\documentclass[12pt]{article}
\usepackage[german, french, american]{babel}
\usepackage[latin1]{inputenc}
\usepackage[T1]{fontenc}
\usepackage{amsmath, amssymb, amsfonts}
\usepackage{amsthm}
\usepackage{mathtools}
\usepackage{mathrsfs}
\usepackage{enumerate}
\usepackage{graphicx}
\usepackage{color}

\usepackage{thumbpdf}
\usepackage{hyperref}

\newtheorem{theorem}{Theorem}
\newtheorem{lemma}[theorem]{Lemma}
\newtheorem{remark}[theorem]{Remark}
\newtheorem{corollary}[theorem]{Corollary}
\newtheorem{definition}[theorem]{Definition}
\newtheorem{assumption}[theorem]{Assumption}

\newcommand{\RR}{\mathbb{R}}
\newcommand{\NN}{\mathbb{N}}

\renewcommand{\Re}{\operatorname{Re}}
\renewcommand{\Im}{\operatorname{Im}}

\numberwithin{equation}{section}

\begin{document}

\title{Strong and Mild Extrapolated $L^{2}$-Solutions to the Heat Equation with Constant Delay}

\author{Denys Khusainov\thanks{Department of Cybernetics, Kiev National Taras Shevchenko University, Kiev, Ukraine \hfill \texttt{d.y.khusainov@gmail.com}} \and
Michael Pokojovy\thanks{Department of Mathematics and Statistics, University of Konstanz, Konstanz, Germany \hfill \texttt{michael.pokojovy@uni-konstanz.de}} \and
Reinhard Racke\thanks{Department of Mathematics and Statistics, University of Konstanz, Konstanz, Germany \hfill \texttt{reinhard.racke@uni-konstanz.de}}}

\date{January 16, 2013}

\maketitle

\begin{abstract}
    We propose a Hilbert space solution theory for a nonhomogeneous heat equation with delay in the highest order derivatives
    with nonhomogeneous Dirichlet boundary conditions in a bounded domain.
    Under rather weak regularity assumptions on the data, we prove a well-posedness result
    and give an explicit representation of solutions.
    Further, we prove an exponential decay rate for the energy in the dissipative case.
    We also show that lower order regularizations lead to ill-posedness, also for higher-order
    equations.
    Finally, an application with physically relevant constants is given.
\end{abstract}

{\bf Key words: } heat equation, delay in highest order terms, strong solutions, mild solutions, well-posedness, ill-posedness

{\bf AMS: }
	35B30, 
	35B35, 
	35D30, 
	35D35, 
	35K20, 
	35Q79  

\pagestyle{myheadings}
\thispagestyle{plain}
\markboth{\uppercase{D. Khusainov, M. Pokojovy, R. Racke}}{\uppercase{Heat Equation with Constant Delay}}

\section{Introduction}
    Let $\Omega \subset \RR^{d}$ be a domain with a Lipschitz-boundary $\partial \Omega$
    and $T > 0$ be a fixed number.
    Let a function $\theta \colon [0, T] \times \bar{\Omega} \to \RR$ denote
    the temperature measured with respect to a reference temperature $\theta_{0}$
    and let $q \colon [0, T] \times \bar{\Omega} \to \RR^{d}$ be the heat flux at a material point $x \in \bar{\Omega}$ at time $t \in [0, T]$.
    With $\rho \colon \bar{\Omega} \to (0, \infty)$ denoting the specific density
    and $c_{\rho} \colon \bar{\Omega} \to (0, \infty)$ denoting the specific heat capacity,
    the energy conservation law reads as
    \begin{equation}
        \rho(x) c_{\rho}(x) \partial_{t} \theta(t, x) + \mathrm{div}\, q(t, x) = h(t, x) \text{ for } t \in (0, T), \; x \in \Omega,
        \label{CONSERVATION_LAW_LINEAR}
    \end{equation}
    where $h$ stands for the intensity of external heat sources.

    To close this equation, a material law postulating a relation between the temperature and the heat flux is required.
    The classical way to do this consists in using Fourier's law of heat conduction stating
    \begin{equation}
        q(t, x) + \lambda(x) \nabla \theta(t, x) = 0 \text{ for } t \in (0, T), \; x \in \Omega, \label{FOURIER_LAW_LINEAR}
    \end{equation}
    where $\lambda \colon \bar{\Omega} \to (0, \infty)$ denotes the heat conductivity being a material property.
    Plugging Equation (\ref{FOURIER_LAW_LINEAR}) into (\ref{CONSERVATION_LAW_LINEAR})
    leads to the classical parabolic heat equation
    \begin{equation}
        \rho(x) c_{\rho}(x) \partial_{t} \theta(t, x)
        - \mathrm{div}\, \big(\lambda(x) \nabla \theta(t, x)\big) = h(t, x) \text{ for } t \in (0, T), \; x \in \Omega.
        \label{LINEAR_HEAT_CONDUCTION_FOURIER}
    \end{equation}
    In many applications, Equation (\ref{LINEAR_HEAT_CONDUCTION_FOURIER}) provides a very accurate macroscopic description
    of the heat conduction phenomenon.
    For some other physical applications,
    the infinite speed of signal propagation arising from Equation (\ref{LINEAR_HEAT_CONDUCTION_FOURIER}) is a significant drawback.

    In particular for these, the following assumption
    \begin{equation}
        q(t, x) + \lambda(x) \nabla \theta(t - \tau, x) = 0 \text{ for } t \in (0, T), \; x \in \Omega \label{DELAY_LAW_LINEAR}
    \end{equation}
    is more realistic from a physical point of view stating that the heat flux notices changes in
    temperature (gradient) not instantaneously, but with some delay.
    The latter leads to the so-called heat equation with pure delay
    \begin{equation}
        \rho(x) c_{\rho}(x) \partial_{t} \theta(t, x)
        - \mathrm{div}\, \big(\lambda(x) \nabla \theta(t - \tau, x)\big) = h(t, x) \text{ for } t \in (0, T), \; x \in \Omega.
        \label{LINEAR_HEAT_CONDUCTION_DELAY}
    \end{equation}
    In addition to severe problems caused by the loss of regularity,
    Equation (\ref{LINEAR_HEAT_CONDUCTION_DELAY}) turns out to be ill-posed (cf. \cite{DreQuiRa2009}).
    One way to overcome this problem is to ``equivalently'' rewrite Equation (\ref{DELAY_LAW_LINEAR}) as
    \begin{equation}
        q(t + \tau, x) + \lambda(x) \nabla \theta(t, x) = 0 \text{ for } t \in (0, T), \; x \in \Omega
        \notag
    \end{equation}
    and perform a {\em formal} Taylor expansion of order one with respect to $\tau$ (cp. \cite{Ca1958}), i.e.,
    \begin{equation}
        \tau q_{t}(t + \tau, x) + q(t, x) + \lambda(x) \nabla \theta(t, x) = 0 \text{ for } t \in (0, T), \; x \in \Omega, \label{CATTANEO_LAW_LINEAR}
    \end{equation}
    to finally obtain
    \begin{equation}
        \begin{split}
            \rho(x) c_{\rho}(x) \partial_{t} \theta(t, x) + \mathrm{div}\, q(t, x) &= h(t, x) \text{ for } t \in (0, T), \; x \in \Omega, \\
            \tau \partial_{t} q(t, x) + q(t, x) + \lambda(x) \nabla \theta(t, x) &= 0 \text{ for } t \in (0, T), \; x \in \Omega.
        \end{split}
        \label{LINEAR_HEAT_CONDUCTION_CATTANEO}
    \end{equation}

    In the present paper,
    we propose another approach to regularize Equation (\ref{LINEAR_HEAT_CONDUCTION_DELAY}).
    For a small parameter $\varepsilon > 0$, we replace Equation (\ref{DELAY_LAW_LINEAR}) with
    \begin{equation}
        q(t, x) + \varepsilon \lambda(x) \nabla \theta(t, x) + \lambda(x) \nabla \theta(t - \tau, x) = 0 \text{ for } t \in (0, T), \; x \in \Omega \label{DELAY_LAW_LINEAR_REGULARIZED}
    \end{equation}
    and arrive at a regularized heat equation
    \begin{equation}
	\begin{split}
		\rho(x) c_{\rho}(x) \partial_{t} \theta(t, x)
        	- \varepsilon \, \mathrm{div}\, \big(\lambda(x) \nabla \theta(t, x)\big) - \mathrm{div}\, \big(\lambda(x) \nabla \theta(t - \tau, x)\big) = h(t, x) \\
		\hfill \text{ for } t \in (0, T), \; x \in \Omega.
	\end{split}
	\label{LINEAR_HEAT_CONDUCTION_DELAY_REGULARIZED}
    \end{equation}
    Though Equation (\ref{LINEAR_HEAT_CONDUCTION_DELAY_REGULARIZED}) is much better behaved than Equation (\ref{LINEAR_HEAT_CONDUCTION_DELAY}),
    standard results on semigroups for delay equations (see, e.g., \cite{BaPia2005}, \cite{Wu1996}) still cannot be applied
    since the delay term is no low order perturbation of the term without delay.
    A semigroup treatment of this problem nevertheless turned out to be possible.
    In \cite{BaSchn2004}, a perturbation result due to Weiss \& Staffans was used to obtain the well-posedness results
    for an even bigger class of equations given by
    \begin{equation}
        \partial_{t} u(t) = \mathcal{A} u(t) + \int_{-r}^{0} \mathrm{d} \mathcal{B}(\theta) u(t + \theta) \text{ for } t > 0, \quad
        u(t) = \varphi(t) \text{ for } t \in [-r, 0], \notag
    \end{equation}
    where $\mathcal{A}$ is a sectorial operator on a Banach space $X$
    and $\mathcal{B} \in \mathrm{BV}\big([-r, 0], \mathcal{L}(D(\mathcal{A}), X)\big)$ has no mass at $0$.
    The particular situation $\mathcal{B} = \eta \mathcal{A} \delta_{-r}$, $\eta \in \mathbb{R}$, was given some additional attention.

    The first systematic treatment of this topic for the case of unbounded operators though probably dates back to \cite{TraWe1976}.
    The authors considered the following evolution equation
    \begin{equation}
        \partial_{t} u(t) = \mathcal{A} u(t) + \mathcal{F}(u_{t}) \text{ for } t > 0, \quad
        u_{0} = \varphi, \notag
    \end{equation}
    where $\mathcal{A}$ is an infinitesimal generator of a $\mathcal{C}^{0}$-semigroup $(S(t))_{t \geq 0}$ on a Banach space $X$,
    $\mathcal{F}$ is a (possibly) unbounded linear or nonlinear operator and $u_{t} = u|_{[t - r, t]}(\cdot - t)$ denotes the history variable.
    In particular, it was shown for the case of $\mathcal{F}$ being a linear differential operator and containing terms of the same order as $\mathcal{A}$
    that the problem possesses a unique mild solution, i.e., a function $u \in H^{1}\big((0, T), X\big)$, $T > 0$, satisfying the integral equation
	\begin{equation}
		u(t) = S(t) \varphi(0) + \int_{0}^{t} S(t - s) \mathcal{F}(u_{s}) \mathrm{d} s \text{ for a.e. } t \in (0, T).
		\notag
	\end{equation}

    In \cite{DiBlKuSi1983}, a similar problem was studied in the strong case, i.e., $u \in H^{1}\big((0, T), X) \cap L^{2}\big((0, T), D(\mathcal{A})\big)$.
    Namely, the authors considered an abstract linear delay equation of the form
    \begin{equation}
        \partial_{t} u(t) = \mathcal{A} u(t) + \mathcal{B} u(t) + \mathcal{L}_{1} u(t - r) + \mathcal{L}_{2} u_{t},
        \text{ for } t > 0 \notag
    \end{equation}
    where $\mathcal{A}$ is a generator of an analytic semigroup on a Hilbert space $H$.
    A typical example of such equation is given by
    \begin{equation}
        \begin{split}
            \partial_{t} u(t, x) &= \partial_{xx} u(t, x) + \partial_{x} u(t, x) +
            \partial_{xx} u(t - r, x) + \int_{-r}^{0} a(s) \partial_{xx} u(t + s, x) \mathrm{d}s \\
		&\hspace{6.5cm} \text{ for } (t, x) \in (0, \infty) \times (0, 1), \\
            u(t, x) &= \varphi(t, x) \text{ for } (t, x) \in (-r, 0) \times (0, 1),
            \quad u(t, 0) = u(t, 1) = 0 \text{ for } t > 0.
        \end{split}
        \notag
    \end{equation}
    Under certain assumptions on the operators $\mathcal{B}$, $\mathcal{L}_{1}$, $\mathcal{L}_{2}$,
    the well-posedness followed from the existence of a semigroup
    associated with the flow $t \mapsto (u(t), u_{t})$.

    In \cite{DiBlKuSi1984}, the authors elaborated on these results by carefully studying
    the $L^{2}$-regularity of the corresponding solution in certain weighted and interpolation spaces and
    presenting a characterization of the infinitesimal generator.

    An $L^{p}$-treatment of delay differential equations with unbounded operators acting on delay terms for $p \in [1, \infty)$ was given in \cite{DiBl2003}.
    In particular, a well-posedness result was obtained for the following problem
    \begin{equation}
        \partial_{t} u(t) = \mathcal{A} u(t) + \mathcal{L} u(t - \tau) \text{ for } t > 0, \notag
    \end{equation}
    where $\mathcal{A}$ is an elliptic operator of order $2m$ and $\mathcal{L}$ is an integro-differential operator of the same order.

    Recently, hyperbolic partial differential equations have also gained a lot of attention.
    In \cite{NiPi2006}, a wave equation with an internal feedback incorporating a delay in the velocity field was studied.
    The initial boundary value problem
    \begin{equation}
        \begin{split}
            \partial_{tt} u(t, x) - \triangle u(t, x) + a_{0} \partial_{t} u(t, x) + a \partial_{t} u(t - \tau, x) &= 0
            \text{ for } (t, x) \in (0, \infty) \times \Omega, \\
            u(t, x) &= 0 \text{ for } (t, x) \in (0, \infty) \times \Gamma_{0}, \\
            \frac{\partial u}{\partial \nu}(t, x) &= 0 \text{ for } (t, x) \in (0, \infty) \times \Gamma_{1}
        \end{split}
        \notag
    \end{equation}
    subject to appropriate initial conditions, where $\Gamma_{0}, \Gamma_{1} \subset \partial \Omega$ are relatively open with $\bar{\Gamma}_{0} \cap \bar{\Gamma}_{1} = \emptyset$,
    was shown to possess a unique strong solution,
    which is exponentially stable if $a_{0} > a > 0$ or instable, otherwise.
    Similar results have also been obtained for the case of a boundary delay.
    This stability study was later carried out in the case of time-varying internal or boundary delay, i.e., $\tau = \tau(t)$,
    in \cite{NiPi2008}, \cite{NiPiVa2011}, etc.

    To the authors' best knowledge, no well-posedness results are available for the case of the delay in Laplacian
    for higher-order in time systems.
    At the same time, replacing stabilizing feedbacks by their delayed counterparts
    are sometimes known to even lead to ill-posedness of the resulting system
    as shown in \cite{Da1997} for the wave-equation and the Euler \& Bernoulli beam.
    The same holds for a general $m$-th order equation with the pure delay (cf. \cite{DreQuiRa2009})
    \begin{equation}
        \partial_{t}^{m} u(t) + \mathcal{A} u(t - \tau) = 0 \text{ for } t > 0 \notag
    \end{equation}
    for an arbitrary unbounded operator $\mathcal{A}$ possessing a sequence of
    eigenvalues $\lambda_{n} \to \infty$, $n \to \infty$,
    or the hyperbolic-parabolic thermoelasticity with pure delay in the second order elliptic part (s. \cite{Ra2012})
    \begin{equation}
        \begin{split}
            a \partial_{tt} u(t, x) - d \partial_{xx} u(t - \tau_{1}, x) + \beta \partial_{x} \theta(t, x) &= 0 \text{ for } (t, x) \in (0, \infty) \times (0, L), \\
            b \partial_{t} \theta(t, x) - k \partial_{xx} \theta(t - \tau_{2}, x) + \beta \partial_{tx} u(t, x) &= 0 \text{ for } (t, x) \in (0, \infty) \times (0, L), \\
            u(t, 0) = u(t, L) = \partial_{x} \theta(t, 0) = \partial_{x} \theta(t, L) &= 0 \text{ for } (t, x) \in (0, \infty)
            \times (0, L).
        \end{split}
        \notag
    \end{equation}
         In the following,
    we propose a natural solution approach in Hilbert spaces
    which employs a generalization of the classical step method for ordinary delay equations (cf. \cite{Go1992})
    rather than the delay semigroup theory.
    In addition to its simplicity and constructivity,
    our approach allows for nonhomogeneous boundary conditions under rather weak regularity assumptions on the boundary data.
    The latter is very useful for various applications in control theory (cf. \cite{NaTe1980}, \cite{Na1986}, \cite{NaYa1989}).
    We want also to point out that our theory can also be applied
    to obtain mild, strong, extrapolated and mild extrapolated solutions in a much more general case
    even in the $L^{p}$-framework with respect to time (cf. Remark \ref{REMARK_SOLUTION_THEORY_GENERAL_CASE}).

    To justify the necessity of the regularization to have at least the same order as the delay term,
    we make essential amendments to the method from \cite{DreQuiRa2009}
    to show that lower order regularizations lead to ill-posedness like in the
    case with pure delay, also for higher-order systems.
    We also refer the reader to \cite{KuSch1983} for a study on necessary conditions for the well-posedness of
    partial differential equations with delay.

    To give an illustration,
    we apply our theory to get a closed form solution
    to a one-dimensional practical problem related to short-pulse laser heating of metal nanofilms
    with physically relevant constants.

\section{Fourier Heat Conduction}
    In this section, we briefly summarize some well-known results for the following
    initial-boundary value problem for the Fourier heat equation with nonhomogeneous Dirichlet boundary conditions
    \begin{equation}
        \begin{split}
		\partial_{t} u(t, x) &=
		\partial_{i} \big(a_{ij}(x) \partial_{j} u(t, x)\big) + b_{i}(x) \partial_{i} u(t, x) + c(x) u(t, x) + \\
		&\phantom{=}\;\; f(t, x) \text{ for } (t, x) \in (0, T) \times \Omega, \\
		u(t, x) &= \gamma(t, x) \text{ for } (t, x) \in (0, T) \times \partial \Omega, \\
		u(0, x) &= u^{0}(x) \text{ for } x \in \Omega.
		\end{split}
		\label{EQUATION_HEAT_CONDUCTION}
    \end{equation}
    Recall that $\partial \Omega$ is assumed to be Lipschitzian throughout the paper.
    To treat the problem from Equation (\ref{EQUATION_HEAT_CONDUCTION}),
    a corresponding operator framework needs to be introduced.
    First, we formally define in the sense of distributions the differential operators
    \begin{equation}
        A_{0} := \partial_{i} \big(a_{ij}(\cdot) \partial_{j}\big), \quad
        A_{r} := b_{i}(\cdot) \partial_{i} + c(\cdot). \notag
    \end{equation}
    Here and in the sequel,
    we employ the Einstein's summation convention.
    So, $\partial_{i} (a_{ij}(\cdot) \partial_{j})$
    should be interpreted as $\sum\limits_{i, j = 1}^{d} \partial_{i} (a_{ij}(\cdot) \partial_{j})$, etc.

    Let $X := L^{2}(\Omega)$ be equipped with the standard scalar product.
    We define the operators
    \begin{equation}
        \begin{split}
            \mathcal{A}_{0} \colon D(\mathcal{A}_{0}) \subset X &\to X, \quad u \mapsto A_{0} u, \\
            \mathcal{A}_{r} \colon D(\mathcal{A}_{r}) \subset X &\to X, \quad u \mapsto A_{r} u
        \end{split}
        \notag
    \end{equation}
    with
    \begin{equation}
        D(\mathcal{A}_{0}) = \big\{u \in H^{1}_{0}(\Omega) \,\big|\, A_{0} u \in X\big\}, \quad
        D(\mathcal{A}_{r}) := H^{1}_{0}(\Omega). \notag
    \end{equation}
    According to \cite[Theorems 9.18 and 12.40]{ReRo2004}, the following assertion holds true.
    \begin{theorem}
        \label{THEOREM_SEMIGROUP}
        Let $\partial \Omega \neq \varnothing$ and let
        $a_{ij} \in W^{1, \infty}(\Omega)$ and $a_{ij} = a_{ji}$, $b_{i}, c_{j} \in L^{\infty}(\Omega)$.
        Further, there may exist a constant $\kappa > 0$ such that
        \begin{equation}
            \operatorname*{ess\,inf}_{x \in \Omega} \xi_{i} a_{ij}(x) \overline{\xi_{j}} \geq \kappa |\xi|^{2} \text{ for all } \xi \in \mathbb{C}^{n}. \notag
        \end{equation}
        Then, the perturbed operator $\mathcal{A} := \mathcal{A}_{0} + \mathcal{A}_{r} \colon D(\mathcal{A}_{0}) \subset X \to X$
        is an infinitesimal generator of an analytic semigroup
        $(S(t))_{t \geq 0}$ on $X$.
    \end{theorem}

    Following the approach described by Lasiecka and Triggiani \cite[Section 0.3]{LaTri2010}
    and taking into account the fact that $\mathcal{A}_{0}$ is continuously invertible,
    we define the extrapolation space $X_{-1}$ as the completion of $X$ with respect to the
    $\|\cdot\|_{-1} := \|\mathcal{A}_{0}^{-1} \cdot\|_{X}$ norm.
    Since $X$ is Hilbertian und therefore reflexive,
    $X_{-1}$ is isomorphic to $(D(\mathcal{A}_{0}^{\ast}))'$.
    Note that $X_{-1}$ is a distributional space, e.g., $X_{-1} \subset H^{-1}(\Omega)$.
    Further, the operator $\mathcal{A}$ can be extended to an operator $\mathcal{A}_{-1} \in L(X, X_{-1})$
    being a generator of an analytic semigroup $(S_{-1}(t))_{t \geq 0}$ of bounded linear operators on $X_{-1}$
    which in its turn is an extension of the semigroup $(S(t))_{t \geq 0}$ from Theorem \ref{THEOREM_SEMIGROUP} onto $X_{-1}$.

    Similar to \cite[Section 3.1]{LaTri2010}, we define the Dirichlet map
    $D \colon L^{2}(\partial \Omega) \to X_{-1}$
    sending each $\gamma \in L^{2}(\partial \Omega)$ to a solution $u \in X_{-1}$ of the problem
    \begin{equation}
        A u = 0 \text{ in } \Omega, \quad u = \gamma \text{ on } \partial \Omega.
        \label{DIRICHLET_MAP}
    \end{equation}

    \begin{lemma}
        \label{LEMMA_BOUNDARY_REGULARITY}
        There holds
        \begin{equation}
            D \in L\big(L^{2}(\partial \Omega), H^{1/2}(\Omega)\big)
            \hookrightarrow L\big(L^{2}(\partial \Omega), X\big)
            \hookrightarrow L\big(L^{2}(\partial \Omega), X_{-1}\big). \notag
        \end{equation}
    \end{lemma}
    \begin{proof}
        See \cite{Gri1972} and \cite{JeKe1995}.
    \end{proof}

    The notion of strong solution from \cite{DiBl2003} in the case of homogeneous boundary data
    motivates the following
    \begin{definition}
        A function $u \in H^{1}\big((0, T), L^{2}(\Omega)\big) \cap L^{2}\big((0, T), H^{2}(\Omega)\big)$
        satisfying Equation (\ref{EQUATION_HEAT_CONDUCTION})
        for a.e. $t \in [0, T]$ is called a strong solution.
    \end{definition}

    \begin{remark}
        \label{REMARK_INITIAL_CONDITION}
        The initial and boundary conditions are satisfied in terms of the continuity of the map
        $u \mapsto \big(u(t^{\ast}), u|_{(0, T) \times \partial \Omega}\big)$,
        \begin{equation}
            \begin{split}
                &H^{1}\big((0, T), L^{2}(\Omega)\big) \cap L^{2}\big((0, T), H^{2}(\Omega)\big) \to \\
                &\hspace{3cm}H^{1}(\Omega) \times \Big(H^{3/4}\big((0, T), L^{2}(\partial \Omega)\big) \cap
                L^{2}\big((0, T), H^{3/2}(\partial \Omega)\big)\Big)
            \end{split}
            \notag
        \end{equation}
        for an arbitrary $t^{\ast} \in [0, T]$ (cf. \cite {Pr2002}).
    \end{remark}

    The fact that a strong solution to Equation (\ref{EQUATION_HEAT_CONDUCTION}) has to satisfy the equation
    \begin{equation}
        \partial_{t} u = \mathcal{A} \big(u - D \gamma\big) + f \text{ in } L^{2}\big((0, T), X\big) \notag
    \end{equation}
    and thus
    \begin{equation}
        \partial_{t} u = \mathcal{A}_{-1} u - \mathcal{A}_{-1} D \gamma + f \text{ in } L^{2}\big((0, T), X_{-1}\big) \notag
    \end{equation}
    motivates the following definition of extrapolated solutions
    (cp. the notion of extrapolated solution in \cite{DiBl2003}).
    \begin{definition}
        \label{DEFINITION_WEAK_SOLUTION}
        A function $u \in H^{1}\big((0, T), X_{-1}\big)$ given by
            \begin{equation}
            \begin{split}
                u(t) = S_{-1}(t) u^{0} - \int_{0}^{t} S_{-1}(t - s) \mathcal{A}_{-1} D \gamma(s) \mathrm{d}s
                + \int_{0}^{t} S_{-1}(t - s) f(s) \mathrm{d}s \\
		\hfill \text{ for a.e. } t \in [0, T]
            \end{split}
            \notag
        \end{equation}
        is called a mild extrapolated solution to Equation (\ref{EQUATION_HEAT_CONDUCTION}).
        If it additionally satisfies $u \in L^{2}\big((0, T), X\big)$,
        we call $u$ a strong extrapolated solution.
    \end{definition}

    \begin{theorem}
        \label{THEOREM_CLASSICAL_PARABOLIC}
        Under the conditions of Theorem \ref{THEOREM_SEMIGROUP},
        Equation (\ref{EQUATION_HEAT_CONDUCTION}) possesses a unique mild extrapolated solution
        if $u_{0} \in X_{-1}$, $f \in L^{2}\big((0, T), X_{-1}\big)$ and $\gamma \in L^{2}\big((0, T), L^{2}(\partial \Omega)\big)$.
        Moreover, if $u_{0} \in X$, $f \in L^{2}\big((0, T), X\big)$ and $\partial \Omega \in \mathcal{C}^{0, 1}$,
        then $u$ is strong extrapolated solution, which additionally satisfies
        \begin{equation}
            u \in L^{2}\big((0, T), H^{1/2}(\Omega)\big) \cap H^{1/4}\big((0, T), L^{2}(\Omega)\big) \cap H^{1}\big((0, T), X_{-1}\big). \notag
        \end{equation}
    \end{theorem}

    \begin{proof}
        See \cite[Section 3.1]{LaTri2010}.
    \end{proof}

    Assuming that $\partial \Omega \in \mathcal{C}^{1, 1}$ and
    exploiting the maximum $L^{p}$-regularity of $\mathcal{A}$ for $p = 2$,
    the following existence and uniqueness theorem follows from \cite{DeHiePr2007}, \cite{Pr2002}.
    In this case, the mild extrapolated solution $u$ is even a strong solution
    and therefore satisfies Equation (\ref{EQUATION_HEAT_CONDUCTION}) pointwise for a.e. $t \in [0, T]$.
    \begin{theorem}
        \label{THEOREM_CLASSICAL_PARABOLIC_HIGHER_REGULARITY}
        Under the conditions of Theorem \ref{THEOREM_SEMIGROUP} and the regularity assumptions
        \begin{equation}
            u^{0} \in H^{1}(\Omega), \;
            f \in L^{2}\big((0, T), L^{2}(\Omega)\big), \;
            \gamma \in H^{3/4}\big((0, T), L^{2}(\partial \Omega)\big) \cap
            L^{2}\big((0, T), H^{3/2}(\partial \Omega)\big) \notag
        \end{equation}
        as well as the compatibility condition
        \begin{equation}
            \gamma(0, \cdot) = u^{0}|_{\partial \Omega}, \notag
        \end{equation}
        the mild extrapolated solution $u$ is a strong solution.
        Moreover, the mapping $(u^{0}, f, \gamma) \mapsto u$
        is an isomorphism between the data space equipped with the corresponding product norm as well as incorporating the compatibility condition
        and the solution space.
    \end{theorem}

    \begin{remark}
        For homogeneous boundary conditions,
        a strong solution
        \begin{equation}
            u \in H^{1}\big((0, T), X\big) \cap L^{2}\big((0, T), D(\mathcal{A})\big) \notag
        \end{equation}
        in sense of \cite{DiBl2001} can be obtained without any extra regularity assumptions on $\partial \Omega$.
        The data have to satisfy
        \begin{equation}
            u^{0} \in \big(X, D(\mathcal{A})\big)_{1/2, 2}, \quad
            f \in L^{2}\big((0, T), X\big), \notag
        \end{equation}
        where the parentheses denote the real interpolation functor.
    \end{remark}

\subsection{Explicit Representation of Solutions}
    \label{SECTION_EXPLICIT_REPRESENTATION_OF_SOLUTIONS}
    In this section, we will briefly outline an explicit solution representation formula for the problem
    \begin{equation}
        \begin{split}
            \partial_{t} u(t) &= \mathcal{A}_{-1} u(t) - \mathcal{A}_{-1} D \gamma(t) + f(t) \text{ for } t \in (0, T), \\
            u(0) &= u^{0} \notag
        \end{split} \label{EQUATION_CLASSICAL_PARABOLIC}
    \end{equation}
    for the case that $\mathcal{A}$ is self-adjoint (i.e., $b \equiv 0$)
    and the data satisfy $u^{0} \in X_{-1}$, $f \in L^{2}\big((0, T), X_{-1}\big)$, $\gamma \in L^{2}\big((0, T), L^{2}(\Omega)\big)$.

    By the virtue of Theorem \ref{THEOREM_CLASSICAL_PARABOLIC},
    the problem possesses a unique mild extrapolated
    solution $u$ given by
    \begin{equation}
        u(t) = S_{-1}(t) u^{0} + \int_{0}^{t} S_{-1}(t - s) (f(s) - D \gamma(s)) \mathrm{d}s. \label{EQUATION_REPRESENTATION_FORMULA_SEMIGROUP}
    \end{equation}
    On the other hand, $\mathcal{A}$ is an elliptic operator having an eigenfunction expansion
    \begin{equation}
        \mathcal{A}u = \sum_{n = 1}^{\infty} \lambda_{n} \langle u, \phi_{n}\rangle_{X} \phi_{n} \text{ for } u \in D(\mathcal{A}), \notag
    \end{equation}
    where $(\lambda_{n})_{n} \subset \RR$, $\lambda_{n} \to -\infty$ for $n \to \infty$
    and $(\phi_{n})_{n} \subset D(\mathcal{A})$ form an orthogonal basis of $X$ (cf. \cite{ReRo2004}).
    Taking into account that the embeddings $D(\mathcal{A}) \hookrightarrow X \hookrightarrow X_{-1}$ are dense and continuous,
    we further obtain
    \begin{equation}
        \mathcal{A}_{-1} u = \sum_{n = 1}^{\infty} \lambda_{n} \langle u, \phi_{n}\rangle_{X_{-1}} \phi_{n} \text{ for } u \in X, \notag
    \end{equation}
    where $\langle \cdot, \cdot\rangle_{X_{-1}}$ is the uniquely defined continuation of $\langle \cdot, \cdot\rangle_{X}$ onto $X_{-1}$.
    Note that $\langle \cdot, \cdot\rangle_{-1}$ and $\langle \cdot, \cdot\rangle_{X_{-1}}$ do not coincide.

    Plugging the ansatz
    \begin{equation}
        u(t) := \sum_{n = 1}^{\infty} u_{n}(t) \phi_{n} \text{ for a.e. } t \in [0, T] \notag
    \end{equation}
    into Equation (\ref{EQUATION_CLASSICAL_PARABOLIC}),
    we obtain a sequence of ordinary differential equations for $u_{n}$
    \begin{equation}
        \begin{split}
            \partial_{t} u_{n}(t) &= \lambda_{n} u_{n}(t) + \langle f(t) - \lambda_{n} D \gamma(t), \phi_{n}\rangle_{X_{-1}} \text{ for a.e. } t \in [0, T], \\
            u_{n}(0) &= \langle u^{0}, u_{n}\rangle_{X_{-1}},
        \end{split}
        \notag
    \end{equation}
    which is uniquely solved by $u_{n} \in H^{1}\big((0, T), \RR\big)$ with
    \begin{equation}
        u_{n}(t) = e^{\lambda_{n} t} \langle u^{0}, \varphi_{n}\rangle_{X_{-1}}
        + \int_{0}^{t} e^{\lambda_{n} (t - s)} \langle f(s) - \lambda_{n} D \gamma(s), \phi_{n}\rangle_{X_{-1}} \mathrm{d}s \text{ for a.e. } t \in [0, T]. \notag
    \end{equation}
    Using Lebesgue's dominated convergence theorem for Bochner integrals,
    we finally obtain for a.e. $t \in [0, T]$
    \begin{equation}
	u(t) = \sum_{n = 1}^{\infty} e^{\lambda_{n} t} \langle u^{0}, \phi_{n}\rangle_{X_{-1}} \phi_{n}
	+ \sum_{n = 1}^{\infty} \int_{0}^{t} e^{\lambda_{n} (t - s)} \langle f(s) - \lambda_{n} D \gamma(s), \phi_{n}\rangle_{X_{-1}} \phi_{n} \mathrm{d}s.
        \label{EQUATION_REPRESENTATION_FORMULA}
    \end{equation}
    Moreover, $u$ coincides with the mild extrapolated solution given in Equation (\ref{EQUATION_REPRESENTATION_FORMULA_SEMIGROUP}).

\subsection{Asymptotical Behavior of Solutions for $t \to \infty$}
    For the sake of completeness, we give a brief discussion on the asymptotics of solutions to Equation (\ref{EQUATION_HEAT_CONDUCTION})
    in the homogeneous case, i.e., $\gamma \equiv 0$, $f \equiv 0$.
    We will be able to generalize these well-known results for the case of regularized heat equation with delay
    in Section \ref{SECTION_ASYMPTOTICAL_BEHAVIOR_DELAY} later on.
    For simplicity, we assume $b_{i} \equiv 0$, $c \equiv 0$ though the exponentially stability easily carries over
    to the case when $\mathcal{A}$ is just positive definite.

    We assume $u^{0} \in X_{-1}$
    and denote by $u_{T}$ for $T > 0$ the mild extrapolated solution to
    \begin{equation}
        \begin{split}
            \partial_{t} u(t, x) &=
            \partial_{i} \big(a_{ij}(x) \partial_{j} u(t, x)\big) \text{ for } (t, x) \in (0, T) \times \Omega, \\
            u(t, x) &= 0 \text{ for } (t, x) \in (0, T) \times \partial \Omega, \\
            u(0, x) &= u^{0}(x) \text{ for } x \in \Omega.
        \end{split}
        \label{EQUATION_HEAT_CONDUCTION_ASYMPTOTICS}
    \end{equation}
    Due to the unique solvability of Equation (\ref{EQUATION_HEAT_CONDUCTION_ASYMPTOTICS}),
    $u_{T_{1}} = u_{T_{2}}|_{[0, T_{1}]}$ for $T_{2} \geq T_{1} > 0$.
    Thus, $(u_{T})_{T > 0}$ can be uniquely continued to a function $u \in \mathcal{C}^{1}\big([0, \infty), X_{-1}\big)$ satisfying
    \begin{equation}
        u(t) = S_{-1}(t) u^{0} \text{ for } t \in [0, \infty). \notag
    \end{equation}

    The energy associated with the solution $u$ is given by
    \begin{equation}
        E(t) := \tfrac{1}{2} \|u(t, \cdot)\|_{X_{-1}}^{2}. \notag
    \end{equation}
    If $u^{0} \in X$, then
    $u(t, \cdot) \in X$ for all $t \geq 0$ (cf. Theorem \ref{THEOREM_CLASSICAL_PARABOLIC_HIGHER_REGULARITY}),
    i.e., $u$ is a classical extrapolated solution (in particular, a strong extrapolated solution), and
    $E(t) = \tfrac{1}{2} \|u(t, \cdot)\|_{X} = \tfrac{1}{2} \int_{\Omega} |u(t, x)|^{2} \mathrm{d}x$
    since $\|\cdot\|_{-1}$ is a continuation of $\|\cdot\|_{X}$.

    \begin{theorem}
        Let $u^{0} \in X_{-1}$. The energy $E$ decays exponentially, i.e.,
        there exists $\omega > 0$ such that
        \begin{equation}
            E(t) \leq e^{-2 \omega t} E(0) \text{ for } t \geq 0. \notag
        \end{equation}
        Moreover, $u \in L^{2}\big((0, \infty), X_{-1}\big)$.
    \end{theorem}

    \begin{proof}
        Using the fact that $(S_{-1}(t))_{t \geq 0}$ is an extension of an exponentially stable semigroup,
        we easily get
        \begin{equation}
             E(t) = \tfrac{1}{2} \|S_{-1} u_{0}\|_{X_{-1}}^{2} \leq
            \tfrac{1}{2} e^{-2\omega t} \|u_{0}\|_{X_{-1}}^{2} = e^{-2\omega t} E(0). \notag
        \end{equation}
        Taking account the measurability of $u$ and estimating
        \begin{equation}
            \int_{0}^{\infty} \|u(t, \cdot)\|_{X_{-1}}^{2} \mathrm{d}t =
            2 \int_{0}^{\infty} E(t) \mathrm{d}t \leq
            2 E(0) \int_{0}^{\infty} e^{-2 \omega t} \mathrm{d}t =
            -\tfrac{E(0)}{\omega} e^{-2 \omega t}\big|_{t = 0}^{t = \infty} = \tfrac{E(0)}{\omega}, \notag
        \end{equation}
        we finally conclude $u \in L^{2}\big((0, \infty), X_{-1}\big)$.
    \end{proof}

\section{Regularized Heat Conduction with Delay}
    \label{SECTION_REGULARIZED_HEAT_CONDUCTION_WITH_DELAY}
    Now, we turn to the heat conduction with constant delay
    \begin{equation}
        \begin{split}
            u_{t}(t, x) &=
            \partial_{i} \big(a_{ij}(x) \partial_{j} u(t, x)\big) + b_{i}(x) \partial_{i} u(t, x) + c(x) u(t, x) + \\
            &\phantom{=}\;\; \partial_{i} \big(\tilde{a}_{ij}(x) \partial_{j} u(t - \tau, x)\big) + \tilde{b}_{i}(x) \partial_{i} u(t - \tau, x) + \tilde{c}(x) u(t - \tau, x) + \\
            &\phantom{=}\;\; f(t, x) \text{ for } (t, x) \in (0, \infty) \times \Omega, \\
            u(t, x) &= \gamma(t, x) \text{ for } (t, x) \in (0, \infty) \times \partial \Omega, \\
            u(0, x) &= u^{0}(x) \text{ for } x \in \Omega, \\
            u(t, x) &= \varphi(t, x) \text{ for } (t, x) \in (-\tau, 0) \times \Omega.
        \end{split}
        \label{EQUATION_HEAT_CONDUCTION_WITH_DELAY}
    \end{equation}

    \begin{assumption}
        We postulate the following conditions.
        \begin{itemize}
            \item Let $\Omega \subset \RR^{d}$ be a bounded with a Lipschitz-boundary.
            \item $a_{ij}, \tilde{a}_{ij} \in W^{1, \infty}(\Omega)$ and $a_{ij} = a_{ji}$, $b_{i}, \tilde{b}_{i}, c, \tilde{c} \in L^{\infty}(\Omega)$.
            \item There exists a constant $\kappa > 0$ such that
            \begin{equation}
                \operatorname*{ess\,inf}_{x \in \Omega} \xi_{i} a_{ij}(x) \overline{\xi_{j}} \geq \kappa |\xi|^{2} \text{ for all } \xi \in \mathbb{C}^{d}. \notag
            \end{equation}
        \end{itemize}
    \end{assumption}

    Similar to the definition of $\mathcal{A}$ in Section 1,
    we define the operator
    \begin{equation}
        \tilde{\mathcal{A}} \colon D(\tilde{\mathcal{A}}) \subset X \to X, \quad
        u \mapsto
        \partial_{i} (\tilde{a}_{ij}(\cdot) \partial_{j} u) + \tilde{b}_{i}(\cdot) \partial_{i} u + \tilde{c}(\cdot) u \notag
    \end{equation}
    with
    \begin{equation}
        D(\tilde{\mathcal{A}}) := \{u \in H^{1}_{0}(\Omega) \,|\, \tilde{\mathcal{A}} u \in X\}. \notag
    \end{equation}
    Further, we need the following assumption:
    \begin{assumption}
        \label{ASSUMPTION_A_A_TILDE}
        Let at least one of the following conditions be fulfilled:
        \begin{itemize}
            \item[{\it i)}] There exists a constant $\tilde \alpha \in \RR \backslash \{0\}$ such that $\tilde{a}_{ij}(x) = \tilde\alpha \, a_{ij}(x)$ for a.e. $x \in \Omega$.
            \item[{\it ii)}] There exists a constant $\tilde{\kappa} > 0$ such that
            \begin{equation}
                \operatorname*{ess\,inf}_{x \in \Omega} \xi_{i} \tilde{a}_{ij}(x) \overline{\xi_{j}} \geq \tilde{\kappa} |\xi|^{2} \text{ for all } \xi \in \mathbb{C}^{d} \notag
            \end{equation}
            and $\partial \Omega$ is of class $\mathcal{C}^{1, 1}$.
        \end{itemize}
    \end{assumption}

    Under Assumption \ref{ASSUMPTION_A_A_TILDE},
    $\tilde{\mathcal{A}}$ is a closed operator.
    Next, we can define $\tilde{X}_{-1} \simeq (D(\tilde{\mathcal{A}}^{\ast}))'$
    and extend $\tilde{\mathcal{A}}$ to an operator $\tilde{\mathcal{A}}_{-1} \in L(X, \tilde{X}_{-1})$.
    Further, one easily gets $D(\mathcal{A}) = D(\tilde{\mathcal{A}})$ and thus $X_{-1} = \tilde{X}_{-1}$.
    Note that due to the elliptic regularity theory the second assumption even implies
    $D(\mathcal{A}) = D(\tilde{\mathcal{A}}) = H^{2}(\Omega) \cap H^{1}_{0}(\Omega)$.

    \begin{remark}
        In fact, the closedness of $\tilde{\mathcal{A}}$ and
        the conditions $D(\mathcal{A}) \subset D(\tilde{\mathcal{A}})$ and $X_{-1} \supset \tilde{X}_{-1}$ would also be sufficient
        for our purposes.
        This is, for example, the case if $\tilde{a}_{ij} \equiv 0$, $\tilde{b}_{i} \equiv 0$
        and $D(\tilde{\mathcal{A}}) = X$.
    \end{remark}

    Following \cite{DiBl2003} for the case of homogeneous data, we introduce the notion of strong solution in the nonhomogeneous case.
    \begin{definition}
        A function $u \in H^{1}\big((0, T), L^{2}\big(\Omega)\big) \cap L^{2}\big((-\tau, T), H^{2}(\Omega)\big)$
        satisfying Equation (\ref{EQUATION_HEAT_CONDUCTION_WITH_DELAY})
        with the boundary and initial conditions interpreted in the sense of Remark \ref{REMARK_INITIAL_CONDITION}
        is called a strong solution.
    \end{definition}

    Likewise, we obtain a formulation of Equation (\ref{EQUATION_HEAT_CONDUCTION_WITH_DELAY}) in the extrapolation space $X_{-1}$
    \begin{equation}
        \begin{split}
            \partial_{t} u(t) &= \mathcal{A}_{-1} u(t) + \tilde{\mathcal{A}}_{-1} u(t - \tau) - \mathcal{A}_{-1} D \gamma(t) + f(t)
            \text{ in } L^{2}\big((0, T), X_{-1}\big), \\
            u(0+) &= u^{0} \text{ in } X_{-1}, \\
            u(t) &= \varphi(t) \text{ in } L^{2}\big((-\tau, 0), X_{-1}\big).
        \end{split}
        \label{EQUATION_HEAT_CONDUCTION_WITH_DELAY_HOM_BD}
    \end{equation}

	\begin{definition}
		A function $u \in L^{2}\big((-\tau, 0), X_{-1}\big) \cap H^{1}\big((0, T), X_{-1}\big)$
		is called a mild extrapolated solution of (\ref{EQUATION_HEAT_CONDUCTION_WITH_DELAY})
		if it satisfies the integro-functional equation
		\begin{equation}
			\begin{split}
				u(t) &= S_{-1}(t) u^{0} + \int_{0}^{t} \left(\tilde{\mathcal{A}}_{-1} S_{-1}(t - s) \big(u(s - \tau) - D \gamma(s)\big)
				+ S_{-1}(t - s) f(s) \right) \mathrm{d}s \\
				&\hspace{3cm} \text{ for a.e. } t \in [0, T], \\
				u(0+) &= u^{0}, \\
				u(t) &= \varphi(t) \text{ for a.e. } t \in [-\tau, 0].
			\end{split}
			\label{EQUATION_HEAT_CONDUCTION_WITH_DELAY_HOM_BD_INTEGRAL}
		\end{equation}
		If $u$ additionally satisfies $u \in L^{2}\big((0, T), X)$,
		we call it a strong extrapolated solution.
	\end{definition}

    We start our considerations by proving the well-posedness in the strong case.
    In contrast to \cite[Theorem 3.3]{DiBlKuSi1984}, no fixed point iteration is required here.
    \begin{theorem}
        \label{THEOREM_EXISTENCE_STRONG_SOLUTION}
        Let $\partial \Omega$ be of class $\mathcal{C}^{1, 1}$. Assume
        \begin{equation}
            \begin{split}
                u^{0} \in H^{1}(\Omega), \quad \varphi \in L^{2}\big((-\tau, 0), H^{2}(\Omega)\big), \\
                \gamma \in H^{3/4}\big((0, T), L^{2}(\partial \Omega)\big) \cap L^{2}\big((0, T), H^{3/2}(\partial \Omega)\big), \quad f \in L^{2}\big((0, T), L^{2}(\Omega)\big)
            \end{split}
            \notag
        \end{equation}
        and $\gamma(0, \cdot) = u^{0}|_{\partial \Omega}$.
        Then the problem (\ref{EQUATION_HEAT_CONDUCTION_WITH_DELAY_HOM_BD}) possesses a unique strong solution $u$.
        Furthermore, there exists a positive constant $C_{T, \tau} > 0$ such that
        \begin{equation}
		\begin{split}
			\|u(t)\|_{X}^{2}
			\leq C_{T, \tau} \Big(\|u^{0}\|_{X}^{2} + \int_{-\tau}^{0} \|\varphi(s)\|_{X}^{2} \mathrm{d}s
			+ \int_{0}^{T} \|\gamma(s)\|_{L^{2}(\partial \Omega)}^{2} \mathrm{d}s
			+ \int_{0}^{T} \|f(s)\|_{X}^{2} \mathrm{d}s\Big) \\
			\hfill \text{ for a.e. } t \in [0, T].
		\end{split}
		\notag
        \end{equation}
    \end{theorem}

    \begin{proof}
        The idea of the proof consists in transforming Equation (\ref{EQUATION_HEAT_CONDUCTION_WITH_DELAY_HOM_BD})
        to an abstract difference equation.
        Without loss of generality, let $T = n\tau$ for some $n \in \NN$.
        Otherwise, consider the problem
        with $f$ and $\gamma$ smoothly continued onto $\left[0, \tau \left\lceil \tfrac{T}{\tau}\right\rceil\right]$.

        We define the operators $r_{k} \colon L^{2}\big((-\tau, T), X_{-1}\big) \to L^{2}\big((0, \tau), X_{-1}\big)$ for $k = 0, \dots, n$ by means of
        \begin{equation}
            (r_{k} g)(s) = u((k - 1) \tau + s) \text{ for } s \in (0, \tau) \notag
        \end{equation}
        for $g \in L^{2}\big((-\tau, T), X_{-1}\big)$ and set
		\begin{equation}
			\begin{split}
				&(u_{0}, \dots, u_{n}) := (r_{0} u, \dots, r_{n} u), \quad
				(f_{1}, \dots, f_{n}) := (r_{1} f, \dots, r_{n} f), \\
				&(\gamma_{1}, \dots, \gamma_{n}) := (r_{1} (D \gamma), \dots, r_{n} (D \gamma)).
			\end{split}
			\notag
		\end{equation}
        Obviously, if $u$ is a strong solution to (\ref{EQUATION_HEAT_CONDUCTION_WITH_DELAY}), then
        $(u_{0}, \dots, u_{n}) \in L^{2}_{\tau}\big((-\tau, T), H^{2}(\Omega)\big)$, $(u_{1}, \dots, u_{n}) \in H^{1}_{\tau}\big((0, T), L^{2}(\Omega)\big) \cap L^{2}_{\tau}\big((0, T), H^{2}(\Omega)\big)$
        by the virtue of Lemma \ref{LEMMA_SEMIDISCRETE_CONTINUOUS_EQUIVALENCE}
        and $(u_{1}, \dots, u_{n})$ solves the following difference-differential equation
        \begin{equation}
            \begin{split}
                \partial_{t} u_{k} &=
                \mathcal{A}_{-1} u_{k} + \tilde{\mathcal{A}}_{-1} u_{k - 1} - \mathcal{A}_{-1} D \gamma_{k} + f_{k}
                \text{ in } L^{2}((0, \tau), X_{-1}) \text{ for } 1 \leq k \leq n, \\
                u_{1}(0) &= u^{0} \text{ in } X, \\
                u_{0} &= \varphi(\cdot + \tau) \text{ in } L^{2}\big((0, \tau), X\big).
            \end{split}
            \label{EQUATION_DIFFERENCE_DIFFERENTIAL}
        \end{equation}
        Next, we show that the converse is also true.
        Let $(u_{0}, \dots, u_{n}) \in L^{2}_{\tau}\big((-\tau, T), H^{2}(\Omega)\big)$ be such that
        $(u_{1}, \dots, u_{n}) \in H^{1}_{\tau}\big((0, T), L^{2}(\Omega)\big)$ solves Equation (\ref{EQUATION_DIFFERENCE_DIFFERENTIAL}).
        According to Lemma \ref{LEMMA_SEMIDISCRETE_CONTINUOUS_EQUIVALENCE} in the Appendix,
        $u \in L^{2}\big((-\tau, 0), H^{2}(\Omega)\big) \cap H^{1}\big((0, T), L^{2}(\Omega)\big)$.
        Exploiting the initial conditions in Equation (\ref{EQUATION_DIFFERENCE_DIFFERENTIAL}), we obtain
        $u(s) = (r_{0} u)(s + \tau) = \varphi(s)$ for a.e. $s \in [-\tau, 0]$
        and $u(0+) = (r_{1} u)(0) = u^{0}$.
        Further, for all $\phi \in \mathcal{C}_{0}^{\infty}\big((0, T), \RR \big)$, we find
	\begin{equation}
		\begin{split}
			\int_{0}^{T} \partial_{t} u(t) \varphi(t) \mathrm{d}t
			&= -\int_{0}^{T} u(t) \partial_{t} \varphi(t) \mathrm{d}t
			= -\sum_{k = 1}^{n} \int_{(k - 1) \tau}^{k \tau} u(t) \partial_{t} \varphi(t) \mathrm{d}t \\
			&= \sum_{k = 1}^{n} \int_{(k - 1) \tau}^{k \tau} \partial_{t} u(t) \varphi(t) \mathrm{d}t
			-\sum_{k = 1}^{n} u(t) \varphi(t)|_{t = (k - 1) \tau}^{t = k \tau} \\
			&= \sum_{k = 1}^{n} \int_{0}^{\tau} \partial_{t} u_{k}(t) \varphi_{k}(t) \mathrm{d}t
			- u(T) \varphi(T) + u(0) \varphi(0)
		\end{split}
		\notag
	\end{equation}
	\begin{equation}
		\begin{split}
			&= \sum_{k = 1}^{n} \int_{0}^{\tau}
			\big(\mathcal{A}_{-1} u_{k}(t) + \tilde{\mathcal{A}}_{-1} u_{k-1}(t) - \mathcal{A}_{-1} D \gamma_{k}(t) + f_{k}(t)\big)
			\varphi_{k}(t) \mathrm{d}t \\
			&= \int_{0}^{T} \big(\mathcal{A}_{-1} u(t) + \tilde{\mathcal{A}}_{-1} u(t - \tau) - \mathcal{A}_{-1} D \gamma(t) + f(t)\big)
			\varphi(t) \mathrm{d}t
		\end{split}
		\notag
	\end{equation}
        with $\varphi_{k} := \varphi((k - 1) \tau + \cdot)$
        meaning that $u$ satisfies Equation (\ref{EQUATION_HEAT_CONDUCTION_WITH_DELAY}) in $X_{-1}$ for a.e. $t \in [0, T]$.
        By the virtue of regularity assumption, $u$ is then a strong solution.
        Now, the existence of a unique strong solution to Equation (\ref{EQUATION_HEAT_CONDUCTION_WITH_DELAY})
        is reduced to the unique solvability of Equation (\ref{EQUATION_DIFFERENCE_DIFFERENTIAL}) in the corresponding space.

        We use mathematical induction to show the latter.
        To start the induction for $k = 1$, we apply Theorem \ref{THEOREM_CLASSICAL_PARABOLIC} to get
        the existence of a unique strong solution $u_{1} \in L^{2}\big((-\tau, 0), H^{2}(\Omega)\big) \cap H^{1}\big((0, \tau), L^{2}(\Omega)\big)$.
        Assume that (\ref{EQUATION_DIFFERENCE_DIFFERENTIAL}) possesses a unique solution
        $(u_{0}, \dots, u_{k}) \in L^{2}_{\tau}\big((-\tau, k\tau), H^{2}(\Omega)\big)$ with
        $(u_{1}, \dots, u_{k}) \in H^{1}_{\tau}\big((0, k\tau), L^{2}(\Omega)\big)$
        for certain $k \in \{1, \dots, n\}$.
        If $k = n$, the claim holds true. Otherwise, $k + 1 \leq n$.
        Taking into account
        \begin{equation}
            H^{1}\big((0, \tau), L^{2}(\Omega)\big) \cap L^{2}\big((0, \tau), H^{2}(\Omega)\big) \hookrightarrow \mathcal{C}^{0}\big([0, \tau], H^{1}(\Omega)\big), \notag
        \end{equation}
        we consider the Cauchy problem
        \begin{equation}
            \begin{split}
                \partial_{t} u_{k+1}(t) &= \mathcal{A}_{-1} u_{k + 1}(t) + g_{k+1}(t) \text{ for } t \in (0, \tau), \\
                u_{k+1}(0) &= u_{k}(\tau) \notag
            \end{split}
        \end{equation}
        with
        \begin{equation}
            g_{k+1} := \tilde{\mathcal{A}}_{-1} u_{k} - \mathcal{A}_{-1} D \gamma_{k+1} + f_{k+1} \in L^{2}\big((0, \tau), L^{2}(\Omega)\big).
        \end{equation}
        By the virtue of Theorem \ref{THEOREM_CLASSICAL_PARABOLIC},
        this problem is uniquely solved by a function $u_{k+1} \in H^{1}\big((0, \tau), L^{2}(\Omega)\big) \cap L^{2}\big((0, \tau), H^{2}(\Omega)\big)$.
        By construction, we have $(u_{0}, \dots, u_{k+1}) \in L^{2}_{\tau}\big((-\tau, (k + 1)\tau), H^{2}(\Omega)\big)$
        and $(u_{1}, \dots, u_{k+1}) \in H^{1}_{\tau}\big((0, (k+1) \tau), L^{2}(\Omega)\big)$.

        There remains to show the a priori estimate.
        From \cite[Theorem 6.3]{Pa1983}, we obtain the existence of
        constants $C_{\sigma} \geq 1$, $C_{\alpha} > 0$ such that
        \begin{equation}
            \|S(t)\|_{L(X, X)} \leq C_{\sigma} \text{ for } t \in [0, \tau]. \notag
        \end{equation}
        Furthermore, by the virtue of \cite[Proposition 0.1]{LaTri2010} there exists a positive constant $C_{\alpha} > 0$ such that
        \begin{equation}
            \int_{0}^{\tau} \|\tilde{\mathcal{A}} S(t) u(t)\|_{X}^{2} \mathrm{d}t \leq
            C_{\alpha} \int_{0}^{\tau} \|u(t)\|_{X}^{2} \mathrm{d}t. \notag
        \end{equation}
        Finally, from Lemma \ref{LEMMA_BOUNDARY_REGULARITY} we get a constant $C_{\gamma} > 0$ such that
        \begin{equation}
            \int_{0}^{\tau} \|D \gamma(t)\|_{X}^{2} \mathrm{d} t \leq
            C_{\gamma} \int_{0}^{\tau} \|\gamma\|_{L^{2}(\partial \Omega)}^{2} \mathrm{d} t. \notag
        \end{equation}

        Applying Duhamel's formula to Equation (\ref{EQUATION_DIFFERENCE_DIFFERENTIAL}),
        we get
        \begin{equation}
            \begin{split}
                u_{1}(t) &= S_{-1}(t) u^{0} +
                \int_{0}^{t} S_{-1}(t - s) \big(\tilde{\mathcal{A}}_{-1} \varphi(s - \tau) + \mathcal{A}_{-1} \gamma_{1}(s) + f_{1}(s)\big) \mathrm{d}s, \\
                u_{k}(t) &= S_{-1}(t) u_{k-1}(\tau) +
                \int_{0}^{t} S_{-1}(t - s) \big(\tilde{\mathcal{A}}_{-1} u_{k-1}(s) + \mathcal{A}_{-1} \gamma_{k}(s) + f_{k}(s)\big) \mathrm{d}s
            \end{split}
            \notag
        \end{equation}
        for a.e. $t \in [0, \tau]$ and $2 \leq k \leq n$ and therefore
        \begin{equation}
            \begin{split}
                \|u_{1}(t)\|_{X}
                &\leq C_{\sigma} \|u^{0}\|_{X} +
                C_{\alpha} \Big(\int_{-\tau}^{0} \|\varphi(s)\|_{X}^{2}\mathrm{d}s\Big)^{1/2} +
                C_{\gamma} \Big(\int_{0}^{\tau} \|\gamma_{1}(s)\|_{L^{2}(\partial \Omega)}^{2}\mathrm{d}s\Big)^{1/2} + \\
                &\phantom{\leq} C_{\sigma} \Big(\int_{0}^{\tau} \|f_{1}(s)\|_{X}^{2} \mathrm{d}s\Big)^{1/2} =: C_{1}, \\
                \|u_{k}(t)\|_{X}
                &\leq
                (C_{\sigma} + C_{\alpha} \sqrt{\tau}) \operatorname*{ess\,sup}_{t \in [0, \tau]} \|u_{k-1}(s)\|_{X} +
                C_{\gamma} \Big(\int_{0}^{\tau} \|\gamma_{k}(s)\|_{L^{2}(\partial \Omega)}^{2} \mathrm{d}s\Big)^{1/2} + \\
                &
                \phantom{\leq}
                C_{\sigma} \Big(\int_{0}^{\tau} \|f_{k}(s)\|_{X}^{2} \mathrm{d}s\Big)^{1/2}
                =: C_{2} \operatorname*{ess\,sup}_{t \in [0, \tau]} \|u_{k-1}(s)\|_{X} + C_{3, k}.
            \end{split}
            \notag
        \end{equation}
        Using discrete Gronwall's lemma (cf. \cite{Ela1999}), we obtain further
        \begin{equation}
            \operatorname*{ess\,sup}_{t \in [0, \tau]} \|u_{k-1}(s)\|_{X} \leq
            \max\{C_{1}, C_{3, k}\} +
            C_{2} \sum_{j = 0}^{k-1} C_{3, k} e^{(k - j - 2) C_{2}}
            \leq
            C_{1} + C_{2} e^{C_{2} k} \sum_{j = 0}^{k} C_{3, k}. \notag
            \notag
        \end{equation}
        Therefore, there exists a constant $C_{T, \tau} > 0$ such that
        \begin{equation}
		\begin{split}
			\|u(t)\|_{X}^{2}
			\leq C_{T, \tau} \Big(\|u^{0}\|_{X}^{2} + \int_{-\tau}^{0} \|\varphi(s)\|_{X}^{2} \mathrm{d}s + \int_{0}^{T} \|\gamma(s)\|_{L^{2}(\partial \Omega)}^{2} \mathrm{d}s
			+ \int_{0}^{T} \|f(s)\|_{X}^{2} \mathrm{d}s\Big) \\
			\hfill \text{ for a.e. } t \in [0, T].
		\end{split}
		\notag
        \end{equation}
        This completes the proof.
    \end{proof}

    \begin{remark}
        Exploiting the isomorphism property from Theorem \ref{THEOREM_CLASSICAL_PARABOLIC_HIGHER_REGULARITY},
        the proof of the previous theorem can be easily amended to obtain the continuous dependence in stronger norms:
        \begin{equation}
		\begin{split}
			\|u\|_{H^{1}((0, T), L^{2}(\Omega)) \cap L^{2}((0, T), H^{2}(\Omega))} \leq C
			\Big(\|u^{0}\|_{H^{1}(\Omega)} + \|\varphi\|_{L^{2}((-\tau, 0), L^{2}(\Omega))} + \\
			\hfill \|f\|_{L^{2}((0, T), L^{2}(\Omega))} +
			\|\gamma\|_{H^{3/4}((0, T), L^{2}(\partial \Omega)) \cap L^{2}((0, T), H^{3/2}(\partial \Omega))}\Big).
		\end{split}
		\notag
        \end{equation}
    \end{remark}

    \begin{remark}
        \label{REMARK_EXISTENCE_STRONG_SOLUTION_HOMOGENEOUS_BC}
        For homogeneous boundary conditions, the proof of Theorem \ref{THEOREM_EXISTENCE_STRONG_SOLUTION}
        can easily be amended to obtain a unique strong solution in sense of \cite{DiBl2003}
        \begin{equation}
            u \in H^{1}\big((0, T), X\big) \cap L^{2}\big((-\tau, T), D(\mathcal{A})\big) \notag
        \end{equation}
        without any additionaly regularity assumptions on $\partial \Omega$ if the data satisfy
        \begin{equation}
            u^{0} \in (X, D(\mathcal{A}))_{1/2, 2}, \quad \varphi \in L^{2}\big((-\tau, 0), D(\mathcal{A})\big), \quad
            f \in L^{2}\big((0, T), X\big). \notag
        \end{equation}
    \end{remark}

    The assumptions of Theorem \ref{THEOREM_EXISTENCE_STRONG_SOLUTION} can be weakened
    if one is interested in strong extrapolated solutions.
    In this case, neither the $\mathcal{C}^{1, 1}$-smoothness of $\partial \Gamma$ nor the compatibilty condition are required.
    Carrying out the proof of Theorem \ref{THEOREM_EXISTENCE_STRONG_SOLUTION} in $X_{-1}$ instead of $X$,
    we get the following result in the extrapolation space.
    \begin{theorem}
            \label{THEOREM_EXISTENCE_STRONG_EXTRAPOLATED_SOLUTION}
        Assume
        \begin{equation}
            u^{0} \in (X, X_{-1})_{1/2, 2}, \; \varphi \in L^{2}\big((-\tau, 0), X\big), \;
            \gamma \in L^{2}\big((0, T), L^{2}(\partial \Omega)\big), \; f \in L^{2}\big((0, T), X_{-1}\big).
            \notag
        \end{equation}
        Then Equation (\ref{EQUATION_HEAT_CONDUCTION_WITH_DELAY_HOM_BD}) possesses a unique strong extrapolated solution $u$.
        Furthermore, there exists a positive constant $C_{T, \tau} > 0$ such that
        \begin{equation}
		\|u(t)\|_{X_{-1}}^{2}
		\leq C_{T, \tau} \Big(\|u^{0}\|_{X_{-1}}^{2} + \int_{-\tau}^{0} \|\varphi(s)\|_{X_{-1}}^{2} \mathrm{d}s +
		\int_{0}^{T} \big(\|\gamma(s)\|_{L^{2}(\partial \Omega)}^{2} + \|f(s)\|_{X_{-1}}^{2}\big) \mathrm{d}s\Big) \notag
        \end{equation}
        for a.e. $t \in [0, T]$.
    \end{theorem}

    Finally, we address the case of mild extrapolated solutions.
    In certain analogy to the proof of Theorem \ref{THEOREM_EXISTENCE_STRONG_SOLUTION},
    we will equivalently transform Equation (\ref{EQUATION_HEAT_CONDUCTION_WITH_DELAY_HOM_BD}) to an integro-difference equation.
    \begin{theorem}
        \label{THEOREM_EXISTENCE_MILD_EXTRAPOLATED_SOLUTION}
        Let
        \begin{equation}
            u^{0} \in X_{-1}, \quad \varphi \in L^{2}\big((-\tau, 0), X_{-1}\big), \quad
            \gamma \in L^{2}\big((0, T), L^{2}(\partial \Omega)\big), \quad f \in L^{2}\big((0, T), X_{-1}\big). \notag
        \end{equation}
        Equation (\ref{EQUATION_HEAT_CONDUCTION_WITH_DELAY_HOM_BD}) possesses a unique mild extrapolated solution $u$.
        Furthermore, there exists a positive constant $C_{T, \tau} > 0$ such that
        \begin{equation}
		\|u(t)\|_{X_{-1}}^{2}
		\leq C_{T, \tau} \Big(\|u^{0}\|_{X_{-1}}^{2} + \int_{-\tau}^{0} \|\varphi(s)\|_{X_{-1}}^{2} \mathrm{d}s +
		\int_{0}^{T} \big(\|\gamma(s)\|_{L^{2}(\partial \Omega)}^{2} + \|f(s)\|_{X_{-1}}^{2}\big) \mathrm{d}s\Big) \notag
        \end{equation}
        for a.e. $t \in [0, T]$.
    \end{theorem}

    \begin{proof}
        Without loss of generality, we assume $T = n\tau$ for a certain $n \in \NN$. Otherwise, consider $f$ and $\gamma$
        trivially continued onto $\big[0, \tau\left\lfloor \frac{T}{\tau}\right\rfloor\big]$.

        With the operators $r_{k}$, $k = 0, \dots, n$, defined in the proof of Theorem \ref{THEOREM_EXISTENCE_STRONG_SOLUTION},
        we let
        \begin{equation}
		\begin{split}
			&(u_{0}, \dots, u_{n}) := (r_{0} u, \dots, r_{n} u), \quad
			(f_{1}, \dots, f_{n}) := (r_{1} f, \dots, r_{n} f), \\
			&(\gamma_{1}, \dots, \gamma_{n}) := (r_{1} (D \gamma), \dots, r_{n} (D \gamma)).
		\end{split}
		\notag
        \end{equation}
        If $u \in L^{2}\big((-\tau, 0), X_{-1}\big) \cap H^{1}\big((0, T), X_{-1}\big)$ is a mild extrapolated solution to Equation (\ref{EQUATION_HEAT_CONDUCTION_WITH_DELAY_HOM_BD}),
        then $(u_{0}, \dots, u_{n}) \in L^{2}_{\tau}\big((-\tau, T), X_{-1}\big)$, $(u_{1}, \dots, u_{n}) \in H^{1}_{\tau}\big((0, T), X_{-1}\big)$
        holds true by the virtue of Lemma \ref{LEMMA_SEMIDISCRETE_CONTINUOUS_EQUIVALENCE}
        and $(u_{1}, \dots, u_{n})$ satisfies the following integro-difference equation
        \begin{equation}
            \begin{split}
                u_{k}(t) &= S_{-1}(t) u_{k}(0) + \int_{0}^{t} \mathcal{A}_{-1} S_{-1}(t - s) \big(u_{k-1}(s) + \gamma_{k}(s)\big) \mathrm{d}s + \\
                &\hspace{2.7cm} \int_{0}^{t} S_{-1}(t - s) f_{k}(s) \mathrm{d}s \text{ for a.e. } t \in [0, \tau], \; 1 \leq k \leq n, \\
                u_{k}(\tau) &= u_{k+1}(0) \text{ for } 1 \leq k \leq n - 1, \\
                u_{1}(0) &= u^{0}, \\
                u_{0} &= \varphi(\cdot + \tau).
            \end{split}
            \label{EQUATION_DIFFERENCE_INTEGRAL}
        \end{equation}
        We claim that the converse is also true.
        Indeed, let $(u_{0}, \dots, u_{n}) \in L^{2}_{\tau}\big((-\tau, T), X_{-1}\big)$ such that
        $(u_{1}, \dots, u_{n}) \in H^{1}_{\tau}\big((0, T), X_{-1}\big)$ solves Equation (\ref{EQUATION_DIFFERENCE_INTEGRAL}).
        Using once again Lemma \ref{LEMMA_SEMIDISCRETE_CONTINUOUS_EQUIVALENCE}, we conclude
        $u \in L^{2}\big((-\tau, T), X_{-1}\big) \cap H^{1}\big((0, T), X_{-1}\big)$.
        From the Equation (\ref{EQUATION_DIFFERENCE_INTEGRAL}) we further deduce
        $u(s) = (r_{0} u)(s + \tau) = \varphi(s)$ for a.e. $s \in [-\tau, 0]$ and
        $u(0+) = (r_{1} u)(0) = u^{0}$.
        There remains to show that the integral equation in (\ref{EQUATION_HEAT_CONDUCTION_WITH_DELAY_HOM_BD_INTEGRAL}) is satisfied.
        This will be shown using mathematical induction.
        For a.e. $t \in [0, \tau]$, we have
        \begin{equation}
            \begin{split}
                u(t) &= S_{-1}(t) u_{1}(0) + \int_{0}^{t} \mathcal{A}_{-1} S_{-1}(t - s) \big(u_{0}(s) + \gamma_{1}(s)\big) \mathrm{d}s +
                \int_{0}^{t} S_{-1}(t - s) f_{1}(s) \mathrm{d}s \\
                &= S_{-1}(t) u^{0} + \int_{0}^{t} \left(\tilde{\mathcal{A}}_{-1} S_{-1}(t - s) \big(u(s - \tau) - D \gamma(s)\big)
                + S_{-1}(t - s) f(s) \right) \mathrm{d}s.
            \end{split}
            \notag
        \end{equation}
        Assume now that the claim is true on $[0, k \tau]$. If $k = n$, the claim trivially holds.
        Otherwise, $k < n$, and we have for a.e. $t \in [k \tau, (k + 1) \tau]$, $\tilde{t} := t - k \tau$
        \begin{equation}
            \begin{split}
                u(t) &= S_{-1}(\tilde{t}) u_{k+1}(0) + \int_{0}^{\tilde{t}} \Big(\mathcal{A}_{-1} S_{-1}(\tilde{t} - s) \big(u_{k}(s) + \gamma_{k+1}(s)\big) + S_{-1}(t - s) f_{k+1}(s)\Big) \mathrm{d}s \\
                &= S_{-1}(\tilde{t}) u_{k}(\tau) + \int_{0}^{\tilde{t}} \mathcal{A}_{-1} S_{-1}(\tilde{t} - s) \big(u_{k}(s) + \gamma_{k+1}(s)\big) \mathrm{d}s +
                \int_{0}^{\tilde{t}} S_{-1}(t - s) f_{k+1}(s) \mathrm{d}s \\
                &= S_{-1}(\tilde{t}) u(k \tau) + \int_{0}^{\tilde{t}} \mathcal{A}_{-1} S_{-1}(\tilde{t} - s) \big(u((k + 1)\tau + s) + D \gamma(k\tau + s)\big) \mathrm{d}s + \\
                &\hspace{7cm} \int_{0}^{\tilde{t}} S_{-1}(t - s) f((k + 1)\tau + s) \mathrm{d}s \\
                &= S_{-1}(\tilde{t}) \bigg(S_{-1}(k \tau) u^{0} + \int_{0}^{k \tau} \Big(\tilde{\mathcal{A}}_{-1} S_{-1}(k \tau - s) \big(u(s - \tau) - D \gamma(s)\big) + \\
		&\hspace{7cm} S_{-1}(k \tau - s) f(s) \Big) \mathrm{d}s\bigg) + \\
                &\phantom{=}\;\; \int_{0}^{\tilde{t}} \mathcal{A}_{-1} S_{-1}(\tilde{t} - s) \big(u(k\tau + s) + D \gamma((k + 1)\tau + s)\big) \mathrm{d}s + \\
                &\phantom{=}\;\;\int_{0}^{\tilde{t}} S_{-1}(t - s) f((k + 1)\tau + s) \mathrm{d}s \\
                &= S_{-1}(t) u^{0} + \int_{0}^{k \tau} \left(\tilde{\mathcal{A}}_{-1} S_{-1}(t - s) \big(u(s - \tau) - D \gamma(s)\big)
                + S_{-1}(t - s) f(s) \right) \mathrm{d}s + \\
                &\phantom{=}\;\; \int_{k \tau}^{t} \mathcal{A}_{-1} S_{-1}(t - s) \big(u(s - \tau) + D \gamma(s)\big) \mathrm{d}s +
                \int_{0}^{t} S_{-1}(t - s) f(s) \mathrm{d}s \\
                &= S_{-1}(t) u_{0} + \int_{0}^{t} \left(\tilde{\mathcal{A}}_{-1} S_{-1}(t - s) \big(u(s - \tau) - D \gamma(s)\big)
                + S_{-1}(t - s) f(s) \right) \mathrm{d}s.
            \end{split}
            \notag
        \end{equation}
        Thus, we have shown that Equations (\ref{EQUATION_HEAT_CONDUCTION_WITH_DELAY_HOM_BD_INTEGRAL}) and (\ref{EQUATION_DIFFERENCE_INTEGRAL})
        are equivalent.

        Again, we exploit mathematical induction to show that Equation (\ref{EQUATION_DIFFERENCE_INTEGRAL}) possesses a unique solution.
        Restricting Equation (\ref{EQUATION_DIFFERENCE_INTEGRAL}) onto $[0, \tau]$, Theorem \ref{THEOREM_CLASSICAL_PARABOLIC} yields the existence of a unique mild extrapolated solution
        \begin{equation}
            u_{1} \in H^{1}\big((0, \tau), X_{-1}\big) = H^{1}_{\tau}\big((0, \tau), X_{-1}\big). \notag
        \end{equation}
        Further, $(u_{0}, u_{1}) \in L^{2}_{\tau}\big((-\tau, 0), X_{-1}\big)$.

        Assume now (\ref{EQUATION_DIFFERENCE_INTEGRAL}) to possess a unique solution
        $(u_{0}, \dots, u_{k}) \in L^{2}\big((-\tau, k\tau), X_{-1}\big)$ such that $(u_{1}, \dots, u_{k}) \in H^{1}_{\tau}\big((0, k \tau), X_{-1}\big)$.
        Excluding the trivial case $k = n$, we have $k < n$.
        Looking at the Equation (\ref{EQUATION_DIFFERENCE_INTEGRAL}) on the $(k + 1)$-st interval
        and exploiting the condition $u_{k+1}(0) = u_{k}(\tau)$, we get
        \begin{equation}
            u_{k+1}(t) = S_{-1}(t) u_{k}(\tau) + \int_{0}^{t} \mathcal{A}_{-1} S_{-1}(t - s) \big(u_{k}(s) + \gamma_{k}(s)\big) \mathrm{d}s +
            \int_{0}^{t} S_{-1}(t - s) f_{k}(s) \mathrm{d}s. \notag
        \end{equation}
        Using the assumptions and the properties of the semigroup $(S_{-1}(t))_{t \geq 0}$, we obtain a unique solution
        \begin{equation}
            u_{k+1} \in H^{1}\big((0, \tau), X_{-1}\big). \notag
        \end{equation}
        Taking into account the condition $u_{k+1}(0) = u_{k}(\tau)$, we finally conclude
        $(u_{0}, \dots, u_{k}) \in L^{2}_{\tau}\big((-\tau, (k + 1) \tau), X_{-1}\big) \cap H^{1}_{\tau}\big((0, (k + 1) \tau), X_{-1}\big)$.
        Thus, the existence proof is finished.

        The proof of continuous dependence on the initial data is literally the same as in the strong case in Theorem \ref{THEOREM_EXISTENCE_STRONG_SOLUTION}
        carried out in $X_{-1}$ instead of $X$.
    \end{proof}

    \begin{remark}
        \label{REMARK_SOLUTION_THEORY_GENERAL_CASE}
        Though this is not the scope of the present paper,
        we want to point out that our method can be applied to a much more general class of problems
        then parabolic ones.
        For a Banach space $X$ and a number $p \in [1, \infty)$, consider the following general delay equation
        \begin{equation}
            \begin{split}
                \partial_{t} u(t) &= \mathcal{A} u(t) + \mathcal{B} u_{t} + f(t) \text{ for } t > 0, \\
                u(0+) &= u^{0}, \\
                u(t) &= \varphi(t) \text{ for } t \in [-\tau, 0]
            \end{split}
            \notag
        \end{equation}
        where $\mathcal{A}$ is a generator of a $\mathcal{C}^{0}$-semigroup of linear, bounded operators on $X$
        and $\mathcal{B} \in L\big(L^{p}\big((-\tau, 0), X\big), L^{p}\big((-\tau, 0), X\big)\big)$.
        Note that $u_{t}$ denotes the usual history variable given by
        $u_{t} \colon [-\tau, 0] \to X$, $s \mapsto u(t + s)$.
        If $\varphi \in L^{p}\big((-\tau, 0), X)$, $f \in L^{p}\big((0, T), X)$,
        same arguments can be exploited to show the existence of a unique mild solution
        $u \in L^{p}\big((-\tau, T), X) \cap W^{1, p}\big((0, T), X)$
        depending continuously on $f$ and $\varphi$.
        Further, the extrapolation space
        $X_{-1}$ can be defined as a completion of $X$ with respect to
        $\|\cdot\|_{-1} := \|(\mathcal{A} + \beta)^{-1} \cdot \|_{X}$, $\beta > 0$ sufficiently large.
        If $X$ is reflexive, the latter can be shown to be isomorphic to $D(\mathcal{A}^{\ast})'$.
        Thus, a mild extrapolated solution $u \in L^{p}\big((-\tau, T), X_{-1}) \cap W^{1, p}\big((0, T), X_{-1})$ can also be constructed.
        To obtain higher regularity for mild solutions or even strong solutions, more knowledge about the structure of $\mathcal{A}$ and $\mathcal{B}$ is though required.
    \end{remark}

\subsection{Explicit Representation of Solutions}
    \label{SUBSECTION_REPRESENTATION_OF_SOLUTIONS}
    In this section, we present an explicit solution formula for Equation (\ref{EQUATION_HEAT_CONDUCTION_WITH_DELAY}).

    For $a, b \in \RR$, we consider first the following scalar ordinary delay differential equation
    \begin{equation}
        \begin{split}
            \partial_{t} u(t) &= a u(t) + b u(t - \tau) + f(t) \text{ for a.e. } t \in [0, T], \\
            u(0) &= u^{0}, \\
            u(t) &= \varphi(t) \text{ for a.e. } t \in [-\tau, 0].
        \end{split}
        \label{EQUATION_ORDINARY_DELAY_DE}
    \end{equation}
    Following the approach in \cite{KhuIvKo2009},
    we define for a number $b \in \RR$ the delayed exponential function
    $\exp_{\tau}(b, \cdot) \colon \RR \to \RR$ given by
    \begin{equation}
        \exp_{\tau}(b, t) :=
        \begin{cases}
            0, & t < -\tau, \\
            1 + \sum\limits_{k = 1}^{\left\lfloor \tfrac{t}{\tau}\right\rfloor + 1} \frac{(t - (k - 1) \tau)^{k}}{k!} b^{k}, & t \geq -\tau.
        \end{cases}
        \notag
    \end{equation}
    Note that the definition can easily be generalized to the case when $b$
    is a matrix or a bounded linear operator on a Banach space $X$.
    \begin{figure}[h!]
        \centering
        \includegraphics[scale = 0.35]{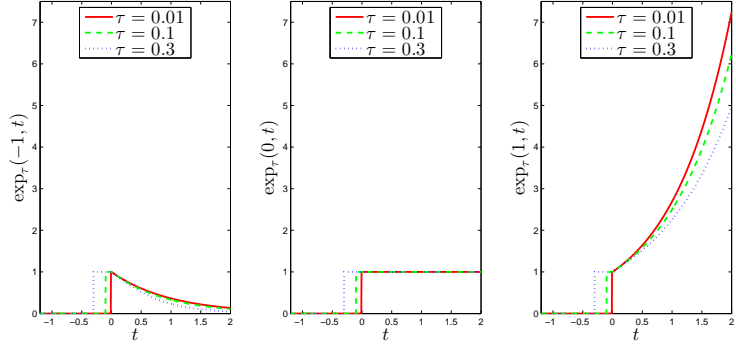}
        \caption{Delayed exponential function}
    \end{figure}
    \begin{theorem}
        \label{THEOREM_DELAY_ODE}
        Let
        $u^{0} \in \RR$, $\varphi \in L^{2}\big((-\tau, 0), \RR)$, $f \in L^{2}\big((0, T), \RR\big)$.
        The delay differential equation (\ref{EQUATION_ORDINARY_DELAY_DE})
        possesses a unique solution
        $u \in L^{2}\big((-\tau, T), \RR\big) \cap H^{1}\big((0, T), \RR\big)$ given by
        \begin{equation}
            \begin{split}
                u(t) &= \left\{
                \begin{array}{cl}
                    \varphi(t), & t \in [-\tau, 0), \\
                    u^{0}, & t = 0, \\
                    {e^{at} \exp_{\tau}(b e^{-a \tau}, t - \tau) u^{0} +
                    b \int_{-\tau}^{0} e^{a(t - s - \tau)} \exp_{\tau}(b e^{-a \tau}, t - 2\tau - s) \varphi(s) \mathrm{d}s + \atop
                    \int_{0}^{t} e^{a(t - s)} \exp_{\tau}(b e^{-a\tau}, t - \tau - s) f(s) \mathrm{d}s,} & t \in (0, T]
                \end{array}\right.
            \end{split}
            \label{EQUATION_SOLUTION_DELAY_ODE}
        \end{equation}
        If $\varphi$ lies in $H^{1}\big((-\tau, 0), \RR)$ and satisfies the compatibility condition $\varphi(0) = u^{0}$,
        $u \in H^{1}\big((-\tau, T), \RR\big)$ holds additionally.
    \end{theorem}
    \begin{proof}
        From \cite{KhuIvKo2009} we know for the classical case,
        i.e., if $\varphi \in \mathcal{C}^{1}\big([-\tau, 0], \RR\big)$ and $\varphi(0) = u^{0}$,
        $f \in \mathcal{C}^{0}\big([-\tau, 0], \RR\big)$,
        that (\ref{EQUATION_ORDINARY_DELAY_DE})
        possesses a unique solution
        $u \in \mathcal{C}^{0}\big([-\tau, T], \RR\big) \cap \mathcal{C}^{1}\big([-\tau, 0], \RR\big) \cap \mathcal{C}^{1}\big([0, T], \RR\big)$ given by
        $u(t) = u_{1}(t) + u_{2}(t)$
        where $u_{1}$ solves (\ref{EQUATION_ORDINARY_DELAY_DE}) for $f \equiv 0$
        and $u_{2}$ solves (\ref{EQUATION_ORDINARY_DELAY_DE}) for $u^{0} = 0$ and $\varphi \equiv 0$.
        It was further shown
        \begin{equation}
            \begin{split}
                u_{1}(t) &=
                \left\{\begin{array}{cl}
                    \varphi(t), & t \in [-\tau, 0), \\
                    u^{0}, & t = 0, \\
                    {\exp_{\tau}(b e^{-a\tau}, t) e^{a(t - \tau)} \varphi(-\tau) + \atop
                    \int_{-\tau}^{0} \exp_{\tau}(b e^{-a\tau}, t - \tau - s) e^{a(t - s)} (\dot{\varphi}(s) - a \varphi(s)) \mathrm{d}s,} &
                    t \in (0, T],
                \end{array}\right. \\
                u_{2}(t) &=
                \left\{\begin{array}{cl}
                    0, & t \in [-\tau, 0], \\
                    \int_{0}^{t} \exp_{\tau}(b e^{-a \tau}, t - \tau - s) e^{a(t - s)} f(s) \mathrm{d}s, & t \in (0, T].
                \end{array}\right.
            \end{split} \notag
        \end{equation}
        Performing partial integration for $\dot{\varphi}$ in $u_{1}$, we obtain for $t \in [0, T]$
	\begin{equation}
		\begin{split}
			u_{1}(t) &= \exp_{\tau}(b e^{-a\tau}, t) e^{a(t + \tau)} \varphi(-\tau) +
			\int_{-\tau}^{0} e^{a(t - s)} \exp_{\tau}(b e^{-a \tau}, t - \tau - s) \dot{\varphi}(s) \mathrm{d}s - \\
			&\phantom{=\;\;} a \int_{-\tau}^{0} e^{a(t - s)} \exp_{\tau}(b e^{-a \tau}, t - \tau - s) \varphi(s) \mathrm{d}s \\
			&= \exp_{\tau}(b e^{-a \tau}, t) e^{a(t + \tau)} \varphi(-\tau) +
			e^{a(t - s)} \exp_{\tau}(b e^{-a \tau}, t - \tau - s) \varphi(s)|_{s = -\tau}^{s = 0} - \\
			&\phantom{=\;\;} \int_{-\tau}^{0} \big(-a e^{a(t - s)} \exp_{\tau}(b e^{-a \tau}, t - \tau - s) - \\
			&\phantom{=\;\;} \phantom{\int_{-\tau}^{0} \big(}
			b e^{a(t - s - \tau)} \exp_{\tau}(b e^{-a \tau}, t - 2\tau - s) \varphi(s)\big) \mathrm{d}s - \\
			&\phantom{=\;\;} a \int_{-\tau}^{0} e^{a(t - s)} \exp_{\tau}(b e^{at}, t - \tau - s) \varphi(s) \mathrm{d}s \\
			&= e^{a(t + \tau)} \exp_{\tau}(b e^{at}, t) \varphi(-\tau) + e^{at} \exp_{\tau}(b e^{at}, t - \tau) \varphi(0) -
			e^{a(t + \tau)} \exp_{\tau}(b e^{at}, t) \varphi(-\tau) - \\
			&\phantom{=\;\;} \int_{-\tau}^{0} \big(-a e^{a(t - s)} \exp_{\tau}(b e^{at}, t - \tau - s) - b e^{a(t - s + \tau)} \exp_{\tau}(b e^{at}, t - 2\tau - s)\big) \varphi(s) -
		\end{split}
		\notag
	\end{equation}
	\begin{equation}
		\begin{split}	
			&\phantom{=\;\;} a \int_{-\tau}^{0} e^{a(t - s)} \exp_{\tau}(b e^{-a \tau}, t - \tau - s) \varphi(s) \mathrm{d}s \\
			&= e^{at} \exp_{\tau}(b e^{-a \tau}, t - \tau) \varphi(0)
			+ b \int_{-\tau}^{0} e^{a(t - s - \tau)} \exp_{\tau}(b e^{-a \tau}, t - 2\tau - s) \varphi(s) \mathrm{d}s.
		\end{split} \notag
        \end{equation}
        Taking now an approximation of $\varphi$ and $f$ with smooth functions,
        we easily deduce the validness of the equation also for the weak case.
    \end{proof}

    To better illustrate Equation (\ref{EQUATION_SOLUTION_DELAY_ODE}),
    we plot solutions to the following scalar delay ordinary differential equation for various values of the parameter $a$:
    \begin{equation}
        \begin{split}
            \partial_{t} u(t) &= a u(t) - u(t - 0.2) + \tfrac{\sin(t)}{1 + t^2} \text{ for } t \in [0, 5], \\
            u(0) &= 1, \\
            u(t) &= e^{-t} \text{ for } t \in [-0.2, 0).
        \end{split} \notag
    \end{equation}
	\vspace{-0.5cm}
    \begin{figure}[h]
        \centering
        \includegraphics[scale = 0.35]{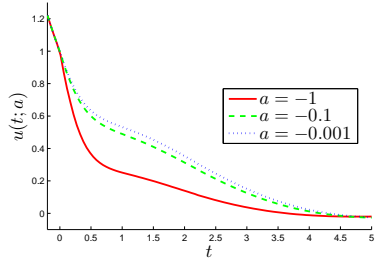}
        \caption{Solution functions $u(\cdot; a)$}
    \end{figure}

    Next, we want to obtain a simple solution representation formula.
    For this purpose, we postulate the following.
    \begin{assumption}
        There exist constants $\alpha \in \RR \backslash \{0\}$, $\beta \in \RR$ such that
        $\tilde{a}_{ij}(x) = \alpha a_{ij}(x)$,
        $\tilde{b}_{i}(x) = b_{i}(x) = 0$ and $\tilde{c}(x) = \alpha c(x) + \beta$ for a.e. $x \in \Omega$.
    \end{assumption}
    Then, both $\mathcal{A}$ and $\tilde{\mathcal{A}}$ are elliptic operators having eigenfunction expansions
    \begin{equation}
        \mathcal{A}u = \sum_{n = 1}^{\infty} \lambda_{n} \langle u, \phi_{n}\rangle_{X} \phi_{n}, \quad
        \tilde{\mathcal{A}} u = \sum_{n = 1}^{\infty} (\alpha \lambda_{n} + \beta) \langle u, \phi_{n}\rangle_{X} \phi_{n}
        \text{ for } u \in D(\mathcal{A}) = D(\tilde{\mathcal{A}}) \notag
    \end{equation}
    with common eigenfunctions forming an orthonormal basis of $X$
    and eigenvalues $(\lambda_{n})_{n} \subset \RR$, $\lambda_{n} \to -\infty$ for $n \to \infty$.
    Similar to Section \ref{SECTION_EXPLICIT_REPRESENTATION_OF_SOLUTIONS},
    we get
    \begin{equation}
        \mathcal{A}_{-1} u = \sum_{n = 1}^{\infty} \lambda_{n} \langle u, \phi_{n}\rangle_{X_{-1}} \phi_{n}, \quad
        \tilde{\mathcal{A}}_{-1} u = \sum_{n = 1}^{\infty} (\alpha \lambda_{n} + \beta) \langle u, \phi_{n}\rangle_{X_{-1}} \phi_{n}
        \text{ for } u \in X. \notag
    \end{equation}

    Plugging the ansatz
    \begin{equation}
        u(t) := \sum_{n = 1}^{\infty} u_{n}(t) \phi_{n} \text{ for a.e. } t \in [-\tau, T] \notag
    \end{equation}
    into Equation (\ref{EQUATION_HEAT_CONDUCTION_WITH_DELAY}),
    we obtain a sequence of ordinary delay differential equations for $u_{n}$
    \begin{equation}
        \begin{split}
            \partial_{t} u_{n}(t) &= \lambda_{n} u_{n}(t) + (\alpha \lambda_{n} + \beta) u_{n}(t - \tau)
            + \langle f(t), \phi_{n}\rangle_{X_{-1}} - \lambda_{n} \langle D \gamma(t), \phi_{n}\rangle_{X_{-1}} \\
            &\hspace{8cm} \text{ for a.e. } t \in [0, T], \\
            u_{n}(0) &= \langle u^{0}, \phi_{n}\rangle_{X_{-1}}, \\
            u_{n}(t) &= \langle \varphi(t), \phi_{n}\rangle_{X_{-1}} \text{ for a.e. } t \in [-\tau, 0].
        \end{split}
        \label{EQUATION_DELAY_ODE_N}
    \end{equation}
    By virtue of Theorem \ref{THEOREM_DELAY_ODE},
    there exists a unique solution $u_{n} \in H^{1}\big((-\tau, T), \RR\big)$ given by
    \begin{equation}
        \begin{split}
            u_{n}(t) &=
            e^{\lambda_{n} t} \exp_{\tau}((\alpha \lambda_{n} + \beta) e^{-\lambda_{n} \tau}, t - \tau) \langle u^{0}, \phi_{n}\rangle_{X_{-1}} + \\
            &\phantom{=\;\;} (\alpha \lambda_{n} + \beta) \int_{-\tau}^{0} e^{\lambda_{n} (t - s - \tau)} \exp_{\tau}((\alpha \lambda_{n} + \beta) e^{-a \tau}, t - 2\tau - s) \langle \varphi(s), \phi_{n}\rangle_{X_{-1}} \mathrm{d}s + \\
            &\phantom{=\;\;} \int_{0}^{t} e^{a(t - s)} \exp_{\tau}(e^{-a\tau}, t - \tau - s) \langle f(s) - \lambda_{n} D \gamma(s), \phi_{n}\rangle_{X_{-1}} \mathrm{d}s
            \text{ for a.e. } t \in [0, T].
        \end{split}
        \label{EQUATION_DELAY_ODE_N_SOL_FORMULA}
    \end{equation}
    Again, using Lebesgue's dominated convergence theorem for Bochner integrals, we find for a.e. $t \in [0, T]$
    \begin{equation}
        \begin{split}
            u(t) &=
            \sum_{n = 1}^{\infty} e^{\lambda_{n} t} \exp_{\tau}((\alpha \lambda_{n} + \beta) e^{-\lambda_{n} \tau}, t - \tau) \langle u^{0}, \phi_{n}\rangle_{X_{-1}} + \\
            &\phantom{=\;\;} \sum_{n = 1}^{\infty} (\alpha \lambda_{n} + \beta) \int_{-\tau}^{0} e^{\lambda_{n} (t - s - \tau)} \exp_{\tau}((\alpha \lambda_{n} + \beta) e^{-a \tau}, t - 2\tau - s) \langle \varphi(s), \phi_{n}\rangle_{X_{-1}} \mathrm{d}s + \\
            &\phantom{=\;\;} \sum_{n = 1}^{\infty} \int_{0}^{t} e^{a(t - s)} \exp_{\tau}(e^{-a\tau}, t - \tau - s) \langle f(s) - \lambda_{n} D \gamma(s), \phi_{n}\rangle_{X_{-1}} \mathrm{d}s.
        \end{split}
        \label{EQUATION_DELAY_ODE_SOL_FORMULA}
    \end{equation}
    Additionally, this function coincides with the mild extrapolated solution
    given in Theorem \ref{THEOREM_EXISTENCE_MILD_EXTRAPOLATED_SOLUTION}.

\subsection{Asymptotical Behavior of Solutions for $t \to \infty$}
    \label{SECTION_ASYMPTOTICAL_BEHAVIOR_DELAY}
    Now we want to study the asymptotics of solutions to Equation (\ref{EQUATION_HEAT_CONDUCTION_WITH_DELAY}).
    For simplicity, we begin our considerations by looking at the case of strong solutions.
    Also, we restrict ourselves to the case $b_{i} = \tilde{b}_{i} \equiv 0$ and
    $c = \tilde{c} \equiv 0$.
    But we point out that a similar result can be obtained if $\mathcal{A}$ and $\tilde{\mathcal{A}}$ are just positive definite.
    With these simplifications, our problem reads as follows
    \begin{equation}
	\begin{split}
		\partial_{t} u(t, x) &=
		\partial_{i} \big(a_{ij}(x) \partial_{j} u(t, x)\big) +
		\partial_{i} \big(\tilde{a}_{ij}(x) \partial_{j} u(t - \tau, x)\big) \\
		&\hspace{3.5cm} \text{ for } (t, x) \in (0, \infty) \times \Omega, \\
		u(t, x) &= 0 \text{ for } (t, x) \in (0, \infty) \times \partial \Omega, \\
		u(0, x) &= u^{0}(x) \text{ for } x \in \Omega, \\
		u(t, x) &= \varphi(t, x) \text{ for } (t, x) \in (-\tau, 0) \times \Omega.
	\end{split}
	\label{EQUATION_HEAT_CONDUCTION_WITH_DELAY_GLOBAL}
    \end{equation}
    We define the energy associated to the solution $u$
    \begin{equation}
        E(t) :=
        \frac{1}{2} \|u(t)\|_{L^{2}(\Omega)}^{2} +
        \frac{1}{2} \int_{t - \tau}^{t} \langle \tilde{a}_{ij}(\cdot) \partial_{j} u(s), \partial_{i} u(s)\rangle_{L^{2}(\Omega)} \mathrm{d}s.
    \end{equation}
	Denote
        \begin{equation}
            \tilde{\lambda} := \|(\tilde{a}_{ij}(\cdot))_{ij}\|_{L^{\infty}(\Omega, \RR^{n \times n})} < \infty. \notag
        \end{equation}

    \begin{theorem}
        \label{THEOREM_EXPONENTIAL_STABILITY}
        Let the Assumption \ref{ASSUMPTION_A_A_TILDE} be satisfied with $\tilde{\alpha} > 0$ and
        let $u^{0} \in (D(A), X)_{1/2, 2}$,
        $\varphi \in L^{2}\big((-\tau, 0), D(A)\big)$.
        There exists then a unique strong solution
        \begin{equation}
            u \in H^{1}_{\mathrm{loc}}\big((0, \infty), X\big) \cap
            L^{2}_{\mathrm{loc}}\big((-\tau, \infty), D(\mathcal{A})\big). \notag
        \end{equation}
        Moreover, if the coefficient matrices $(a_{ij}(\cdot))_{ij}$, $(\tilde{a}_{ij}(\cdot))_{ij}$ are such that the condition
	\begin{equation}
		\kappa > \tilde{\lambda} \sqrt{\tfrac{\tilde{\lambda}}{\tilde{\kappa}}} \notag
	\end{equation}
	is satisfied, there exist constants $\omega, C > 0$, independent from the initial data, with
        \begin{equation}
            E(t) \leq C e^{-2 \omega t} E(0) \text{ for a.e. } t \in [0, \infty). \notag
        \end{equation}
    \end{theorem}

    \begin{proof}
        From Remark \ref{REMARK_EXISTENCE_STRONG_SOLUTION_HOMOGENEOUS_BC} we obtain
        the existence of a strong solution on each finite time inverval.
        The latter can thus be continued to a global strong solution
        \begin{equation}
            u \in H^{1}_{\mathrm{loc}}\big((0, \infty), X\big) \cap
            L^{2}_{\mathrm{loc}}\big((-\tau, \infty), D(\mathcal{A})\big). \notag
        \end{equation}
        Note that no regularity assumptions on $\partial \Omega$ are required here
        since homogeneous boundary conditions are considered.
        The following calculations should be interpreted in $L^{2}_{\mathrm{loc}}\big((0, \infty), X)$.

        Multiplying Equation (\ref{EQUATION_HEAT_CONDUCTION_WITH_DELAY_GLOBAL})
        with $u(t, \cdot)$ in $L^{2}(\Omega)$ and carrying out a partial integration yields
        \begin{equation}
                \frac{1}{2} \partial_{t} \|u(t, \cdot)\|_{L^{2}(\Omega)}^{2} =
                - \langle a_{ij}(\cdot) \partial_{j} u(t, \cdot), \partial_{i} u(t, \cdot)\rangle_{L^{2}(\Omega)}
                - \langle \tilde{a}_{ij}(\cdot) \partial_{j} u(t - \tau, \cdot), \partial_{i} u(t, \cdot)\rangle_{L^{2}(\Omega)}.
            \notag
        \end{equation}
        Taking into account the uniform positive definiteness of $a$,
        we get for an arbitrary number $\varepsilon > 0$
        \begin{equation}
            \begin{split}
                \frac{1}{2} \partial_{t} \|u(t, \cdot)\|_{L^{2}(\Omega)}^{2} \leq
                - \kappa \|\nabla u(t, \cdot)\|_{L^{2}(\Omega)}^{2}
                + \tfrac{\tilde{\lambda} \varepsilon}{2} \| \nabla u(t, \cdot)\|_{L^{2}(\Omega)}^{2} + \\
                \hfill \tfrac{\tilde{\lambda}}{2 \varepsilon} \|\nabla u(t - \tau, \cdot)\|_{L^{2}(\Omega)}^{2}.
            \end{split}
            \label{EQUATION_HEAT_CONDUCTION_WITH_DELAY_GLOBAL_MULTIPLIED_AND_ESTIMATED}
        \end{equation}

        Following the standard approach, we consider the history variable
        \begin{equation}
            z(s, t, \cdot) := u(t - s \tau, \cdot) \text{ for } s \in [0, 1], \; t \in [0, \infty).
            \label{EQUATION_HISTORY_FUNCTION_DEFINITION}
        \end{equation}
        Then, $z$ is smooth in $s$ and $t$ (cp. \cite[Lemma 3.4]{BaPia2005}) and there holds in the distributional sense
        \begin{equation}
            \tau \partial_{t} z(s, t, \cdot) + \partial_{s} z(s, t, \cdot)
             = 0 \text{ for a.e. } (s, t) \in (0, 1) \times (0,
             \infty).
            \label{EQUATION_HISTORY_FUNCTION_EVOLUTION}
        \end{equation}
        Further, a transformation of variables yields
        \begin{equation}
            \int_{t - \tau}^{t} \langle \tilde{a}_{ij}(\cdot) \partial_{j} u(s), \partial_{i} u(s)\rangle_{L^{2}(\Omega)} \mathrm{d}s =
            \tau \int_{0}^{1} \langle \tilde{a}_{ij}(\cdot) \partial_{j} z(s, t, \cdot), \partial_{i} z(s, t, \cdot) \rangle_{L^{2}(\Omega)} \mathrm{d}s.
            \notag
        \end{equation}

        For a smooth nonnegative weight function $\rho \colon [0, \tau] \to \RR$ to be
        selected later,
        we define the functional
        \begin{equation}
            F(t) := \int_{0}^{1} \rho(\tau s)
            \left\langle \tilde{a}_{ij}(x) \partial_{j} z(s,t,\cdot), \partial_{i}
            z(s,t,\cdot) \right\rangle_{L^{2}(\Omega)} \mathrm{d}s.
        \end{equation}
        Exploiting Equation (\ref{EQUATION_HISTORY_FUNCTION_EVOLUTION})
	and the identity
	\begin{equation}
		\begin{split}
			\partial_{s} \Big(\rho(\tau s) \langle &\partial_{i} \big(\tilde{a}_{ij}(x) \partial_{j} z(s, t, \cdot), z(s, t, \cdot)\rangle_{L^{2}(\Omega)}\Big) = \\
			&\tau \rho'(\tau s) \langle \partial_{i} \big(\tilde{a}_{ij}(x) \partial_{j} z(s, t, \cdot), z(s, t, \cdot)\rangle_{L^{2}(\Omega)} + \\
			&2 \Re \rho(\tau s) \langle \partial_{i} \big(\tilde{a}_{ij}(x) \partial_{j} z(s, t, \cdot), \partial_{s} z(s, t, \cdot)\rangle_{L^{2}(\Omega)},
		\end{split}
		\notag
	\end{equation}
	we obtain
	\begin{equation}
		\begin{split}
			\frac{\mathrm{d}}{\mathrm{d}t} F(t) &=
			-\frac{2}{\tau} \Re \int_{0}^{1} \rho(\tau s) \langle \partial_{i} \big(\tilde{a}_{ij}(x) \partial_{j} z(s, t, \cdot), \partial_{s} z(s, t, \cdot)\rangle_{L^{2}(\Omega)} \\
			&= -\int_{0}^{1} \rho'(\tau s) \langle \partial_{i} \big(\tilde{a}_{ij}(x) \partial_{j} z(s, t, \cdot),  z(s, t, \cdot)\rangle_{L^{2}(\Omega)} + \\
			&\phantom{=\;\;} \frac{1}{\tau} \int_{0}^{1} \partial_{s} \Big(\rho(\tau s) \langle \partial_{i} \big(\tilde{a}_{ij}(x) \partial_{j} z(s, t, \cdot), z(s, t, \cdot)\rangle_{L^{2}(\Omega)}\Big) \mathrm{d}s \\
			&= \int_{0}^{1} \rho'(\tau s) \left\langle \tilde{a}_{ij}(x) \partial_{j} z(s, t, \cdot), \partial_{i} z(s, t, \cdot) \right\rangle_{L^{2}(\Omega)} \mathrm{d}s - \\
                	&\phantom{=\;\;} \frac{1}{\tau} \Big(\rho(\tau) \left\langle \tilde{a}_{ij}(x) \partial_{j} z(1,t,\cdot), \partial_{i} z(1,t,\cdot)\right\rangle_{L^{2}(\Omega)} - \\
			&\phantom{=\; \frac{1}{\tau}\Big)}
			\rho(0) \langle \tilde{a}_{ij}(x) \partial_{j} z(0, t, \cdot), \partial_{i} z(0, t, \cdot)\rangle_{L^{2}(\Omega)}\Big).
		\end{split}
		\notag
	\end{equation}
        Assuming that $\rho$ is strictly monotonically decreasing, letting
        \begin{equation}
            \rho_{0} := -\tfrac{1}{\tau} \mathop{\max}_{s \in [0, \tau]} \rho'(s) \notag
        \end{equation}
        and exploiting the uniform positive definiteness of $\tilde{a}$, we obtain the estimate
        \begin{equation}
            \begin{split}
                \frac{\mathrm{d}}{\mathrm{d}t} F(t) &\leq
                -\rho_{0} \int_{t - \tau}^{t} \left\langle \tilde{a}_{ij}(x) \partial_{j} u(s, \cdot), \partial_{i} u(s, \cdot) \right\rangle_{L^{2}(\Omega)} \mathrm{d}s - \\
                &\phantom{=\;-} \tfrac{\tilde{\kappa} \rho(\tau)}{\tau} \|\nabla u(t - \tau, \cdot)\|_{L^{2}(\Omega)}^{2} + \tfrac{\tilde{\lambda} \rho(0)}{\tau}
                \|\nabla u(t, \cdot)\|_{L^{2}(\Omega)}^{2}.
            \end{split}
            \label{EQUATION_LYAPUNOV_FUNCTIONAL}
        \end{equation}
        Now, we can define the Lyapunov functional
        \begin{equation}
            L(t) := \frac{1}{2} \|u(t, \cdot)\|_{L^{2}(\Omega)}^{2} + F(t). \notag
        \end{equation}
        Combining (\ref{EQUATION_HEAT_CONDUCTION_WITH_DELAY_GLOBAL_MULTIPLIED_AND_ESTIMATED}) and
        (\ref{EQUATION_LYAPUNOV_FUNCTIONAL}) we obtain
        \begin{equation}
            \begin{split}
                \frac{\mathrm{d}}{\mathrm{d}t} L(t) &\leq
                -\alpha_{1} \|u(t, \cdot)\|_{L^{2}(\Omega)}^{2}
                -\alpha_{2} \|u(t - \tau, \cdot)\|_{L^{2}(\Omega)}^{2} - \\
		&\phantom{=\;\;-} \rho_{0} \int_{t - \tau}^{t} \left\langle \tilde{a}_{ij}(x) \partial_{j} u(s, \cdot), \partial_{i} u(s, \cdot) \right\rangle_{L^{2}(\Omega)} \mathrm{d}s,
                \notag
            \end{split}
        \end{equation}
        where
        \begin{equation}
            \alpha_{1} := \kappa - \tfrac{\tilde{\lambda} \varepsilon}{2} - \tfrac{\tilde{\lambda} \rho(0)}{\tau}, \quad
            \alpha_{2} := \tfrac{\tilde{\kappa} \rho(\tau)}{\tau} - \tfrac{\tilde{\lambda}}{2 \varepsilon}.
            \notag
        \end{equation}
        Now, we have to select $\varepsilon$
        and a smooth, uniformly positive function $\rho \colon [0, 1] \to \mathbb{R}$, e.g., a linear function being uniquely determined by prescribing $\rho(0)$ and $\rho(\tau)$,
	such that $\rho_{0}, \alpha_{1}, \alpha_{2}$ are positive.
	This yields a system of three inequalities
	\begin{equation}
		\kappa - \tfrac{\tilde{\lambda} \varepsilon}{2} - \tfrac{\tilde{\lambda} \rho(0)}{\tau} > 0, \quad
		\tfrac{\tilde{\kappa} \rho(\tau)}{\tau} - \tfrac{\tilde{\lambda}}{2 \varepsilon} > 0, \quad
		\rho(0) > \rho(\tau).
		\notag
	\end{equation}
	After some simple equivalent transformations, we obtain
	\begin{equation}
		\rho(0) > \rho(\tau) > \tfrac{\tau \tilde{\lambda}}{2 \varepsilon \tilde{\kappa}}, \quad
		\kappa > \tfrac{\tilde{\lambda}}{2} \left(\varepsilon + \tfrac{2 \rho(0)}{\tau}\right)
		\label{EQUATION_CONSTANTS_INEQUALITIES}
	\end{equation}
	and thus
	\begin{equation}
		\kappa > \tfrac{\tilde{\lambda}}{2} \left(\varepsilon + \tfrac{\tilde{\lambda}}{\varepsilon \tilde{\kappa}}\right) =: \chi(\varepsilon). \notag
	\end{equation}
	The function $\chi$ attains its global minimum over $\varepsilon > 0$ in $\varepsilon^{\ast} = \sqrt{\frac{\tilde{\lambda}}{\tilde{\kappa}}}$
	with $\chi(\varepsilon^{\ast}) = \tilde{\lambda} \sqrt{\tfrac{\tilde{\lambda}}{\tilde{\kappa}}}$.
	Plugging now $\varepsilon = \varepsilon^{\ast}$ into Equation (\ref{EQUATION_CONSTANTS_INEQUALITIES}), we finally get the ``optimal'' conditions
	\begin{equation}
		\rho(0) > \rho(\tau) > \tfrac{\tau}{2} \sqrt{\tfrac{\tilde{\lambda}}{\tilde{\kappa}}}, \quad
		\kappa > \tilde{\lambda} \sqrt{\tfrac{\tilde{\lambda}}{\tilde{\kappa}}}. \notag
	\end{equation}
	The first inequality can be satisfied by a proper choice of $\rho(0)$ and $\rho(\tau)$.
	The validness of the second inequality is guaranteed by the assumptions. Thus, we have
	$\beta := \min\{\alpha_{1}, \alpha_{2}, \rho_{0}\} > 0$
        and therefore
        \begin{equation}
            \frac{\mathrm{d}}{\mathrm{d}t} L(t) \leq -\beta E(t).
            \label{EQUATION_DERIVATION_LYAPUNOV_VS_ENERGY_ESTIMATE}
        \end{equation}
        Exploiting the monotonicity of $\rho$, we find
        \begin{equation}
            \min\left\{1, \tfrac{2 \rho(\tau)}{\tau}\right\} E(t) \leq L(t) \leq \max\left\{1, \tfrac{2 \rho(0)}{\tau}\right\} E(t) \text{ for a.e. } t \in [0, \infty).
            \label{EQUATIONLYAPUNOV_ENERGY_EQUIVALENCE}
        \end{equation}
        Combining (\ref{EQUATION_DERIVATION_LYAPUNOV_VS_ENERGY_ESTIMATE})
        and (\ref{EQUATIONLYAPUNOV_ENERGY_EQUIVALENCE}), we further arrive at
        \begin{equation}
            \frac{\mathrm{d}}{\mathrm{d}t} L(t) \leq -2\omega L(t) \text{ for a.e. } t \in [0, \infty)
            \notag
        \end{equation}
        with $\omega := \tfrac{\beta}{2} \min\left\{1, \tfrac{\tau}{2 \rho(0)}\right\}$.
        Gronwall's inequality now yields
        \begin{equation}
            L(t) \leq e^{-2\omega} L(0). \notag
        \end{equation}
        Exploiting once again Equation (\ref{EQUATIONLYAPUNOV_ENERGY_EQUIVALENCE}),
        the claim follows with $\omega$ as above and $C := \tfrac{\rho(0)}{\rho(\tau)}$.
    \end{proof}

    Taking into account the equivalence of
    $u \mapsto \big(\langle \tilde{a}_{ij}(\cdot) u, u\rangle_{L^{2}(\Omega)}\big)^{1/2}$ and the norms of interpolation spaces
    $(X, D(\tilde{\mathcal{A}}))_{1/2, 2}$, $(X, D(\mathcal{A}))_{1/2, 2}$,
    the energy $E$ can easily be seen to be equivalent with the squared norm of
    $X \times L^{2}\big((-\tau, 0), (X, D(\mathcal{A}))_{1/2, 2}\big)$.
    Using the extrapolation methods, the energy can thus be continously extended onto
    $X_{-1} \times L^{2}\big((-\tau, 0), (X_{-1}, X)_{1/2, 2}\big)$.
    By approximating the initial data with regular functions and applying Theorem \ref{THEOREM_EXPONENTIAL_STABILITY}, we get the following
    \begin{corollary}
        Let the Assumption \ref{ASSUMPTION_A_A_TILDE} be satisfied and
        let $u^{0} \in (X_{-1}, X)_{1/2, 2}$,
        $\varphi \in L^{2}\big((-\tau, 0), X\big)$.
        There exists then a unique strong extrapolated solution
        \begin{equation}
            u \in H^{1}_{\mathrm{loc}}\big((0, \infty), X_{-1}\big) \cap
            L^{2}_{\mathrm{loc}}\big((-\tau, \infty), X\big). \notag
        \end{equation}
        Moreover, there exist constants $\omega, C > 0$ independent from the initial data
        such that
        \begin{equation}
		\begin{split}
			\|u(t)\|_{X_{-1}}^{2} + \int_{t - \tau}^{t} \|u(s)\|_{(X_{-1}, X)_{1/2, 2}}^{2} \mathrm{d}s \leq
			C e^{-2\omega t} \Big(\|u^{0}\|_{X_{-1}}^{2} + \\
			\hfill \|\varphi\|_{L^{2}((-\tau, 0), (X_{-1}, X)_{1/2, 2})}^{2}\Big)
			\text{ for a.e. } t \geq 0.
		\end{split}
		\notag
        \end{equation}
    \end{corollary}

\section{Ill-Posedness for Lower Order Regularizations}
    To justify the ``sharpness'' of the results from Section \ref{SECTION_REGULARIZED_HEAT_CONDUCTION_WITH_DELAY},
    we show that lower order regularizations of the heat equation with pure delay (\ref{LINEAR_HEAT_CONDUCTION_DELAY})
    lead to an ill-posed problem.
    \begin{theorem}
        \label{THEOREM_ILL_POSEDNESS_LOW_ORDER_REG_HEAT_EQUATION}
        Let $\mathcal{A}$ be defined as in previous section and let $\alpha \in [0, 1)$, $\varepsilon > 0$.
        Let $u^{0} \in (X, D(\mathcal{A}))_{1/2, 2}$, $\varphi \in L^{2}\big((-\tau, 0), D(\mathcal{A})\big)$.
        Then the problem
        \begin{equation}
            \begin{split}
                \partial_{t} u(t) &= -\varepsilon (-\mathcal{A})^{\alpha} u(t) + \mathcal{A} u(t - \tau) \text{ for } t \in (0, T), \\
                u(0) &= u^{0}, \\
                u(t) &= \varphi(t) \text{ for } t \in (-\tau, 0)
            \end{split}
            \notag
        \end{equation}
        is ill-posed.
    \end{theorem}

    More generally, we prove
    \begin{theorem} \label{theo-m1}
        Let $\mathcal{A}$ be a self-adjoint positive operator having a complete orthonormal set of eigenfunctions
        $(\phi_j)_{j}$ corresponding to eigenvalues $(\tilde\lambda_j)_{}$
        with $\tilde\lambda_j \to -\infty$ as $j \to \infty$. Let $\varepsilon > 0$, and let $\alpha \in (-\infty, 1)$.
        Then, the problem
        \begin{eqnarray*}
            \partial_t u(t) &=& {\cal A} u(t-\tau) - \varepsilon (-{\cal A})^\alpha u(t),\\
            u(0)&=&u_0\in (X, D(\mathcal{A}))_{1/2, 2},\\
            u(t)&=&\varphi(t)\qquad \mbox{for } t\in (-\tau,0)\; \mbox{and }
            \varphi\in L^{2}\big((-\tau, 0), D(\mathcal{A})\big)
        \end{eqnarray*}
        is ill-posed. That is, there exists solutions $(u_j)_j$ with norm
        $\|u_j(t)\|$, $j\in {\Bbb N}$, such that, for any fixed $t>0$, the norm
        tends to infinity as $j\to\infty$, while the norm of the data $(u_j(0), \varphi_j)$ remains bounded.
    \end{theorem}

    \begin{proof}
        We make the ansatz
        \begin{equation}
            \label{eq3.1}
            u_j(t) = e^{\omega_j t} \phi_j,
        \end{equation}
        looking for suitable $\omega_j$ such that $\Re \omega_j\to\infty$.
        For such solutions, the norm of the corresponding data will remain bounded,
        but, for any $t > 0$,
        $\|u_j(t)\| = e^{\Re \omega_j t} \to \infty$ as $j \to \infty$.

        The ansatz (\ref{eq3.1}) yields a solution if
        \begin{equation}
            \label{eq3.2}
            \omega_j = -\lambda_j e^{-\tau \omega_j} - \varepsilon
            \lambda_j^\alpha,
        \end{equation}
        where $\lambda_j:=-\tilde\lambda_j\to +\infty$.
        For simplicity, we drop the index $j$ and define
        \begin{equation}\
            \label{eq3.3}
            v := \omega + \varepsilon \lambda^\alpha.
        \end{equation}
        Then, $v$ should satisfy
        \begin{equation}
            \label{eq3.4}
            v = -\lambda e^{\tau \varepsilon \lambda^\alpha} e^{-\tau v}.
        \end{equation}
        Recalling the proof of Theorem 1.1 from \cite{DreQuiRa2009}, there are solutions to (\ref{eq3.4}) satisfying
        \begin{equation}
           \Re v \to \infty \qquad \mbox{ as } \lambda \to \infty.
            \notag
        \end{equation}
        We shall show that
        \begin{equation}
            \label{eq3.5}
            \Re \omega = \Re v - \varepsilon\lambda^\alpha \to \infty
        \end{equation}
        is also valid. This is obvious if $\alpha \leq 0$, therefore, it remains to consider the case $\alpha \in (0, 1)$.
        Observing
        \begin{equation}
            |v| = \lambda e^{\tau(\varepsilon\lambda^\alpha - \Re v)},
            \notag
        \end{equation}
        we obtain
        \begin{equation}
            \Re v = \varepsilon \lambda^\alpha - \tfrac{1}{\tau} \ln(\tfrac{|v|}{\lambda})
            \notag
        \end{equation}
        and
        \begin{equation}
            \Re \omega = \tfrac{1}{\tau} \ln(\tfrac{\lambda}{|v|}).
            \notag
        \end{equation}
        Hence, (\ref{eq3.5}) is equivalent to proving
        \begin{equation}
            \label{eq3.6}
            \tfrac{|v|}{\lambda} \to 0 \qquad \mbox{ as } \lambda\to \infty.
        \end{equation}
        Since we conclude from \cite{DreQuiRa2009} that
        \begin{equation}
            \Im v \to \tfrac{\pi}{\tau},
            \notag
        \end{equation}
        this is equivalent to proving
        \begin{equation}
            \label{eq3.7}
            \tfrac{\Re v}{\lambda} \to 0.
        \end{equation}
        It is interesting to notice that we shall prove that $\Re v$ goes
        to infinity faster than the power term $\lambda^\alpha$ (cp. (\ref{eq3.5}))
        by proving that $\Re v$ goes less fast to infinity than the power term $\lambda$ (cp. (\ref{eq3.7})).
        This will be, of course, possible only
        because $\alpha < 1$ holds.

        To prove (\ref{eq3.7}), we apply the rule of de l'Hospital to $v$ as a
        function in $\lambda$. The relation (\ref{eq3.4}) implies for the
        derivative of $v$
        \begin{equation}
            v' e^{\tau v}(1 + \tau v) = -e^{\tau
            \varepsilon \lambda^\alpha}(1 + \tau \varepsilon\alpha \lambda^\alpha)
            \notag
        \end{equation}
        or
        \begin{equation}
            \begin{split}
                v' &= -\frac{e^{\tau \varepsilon \lambda^\alpha}(1 + \tau \varepsilon \alpha \lambda^\alpha)}{e^{\tau v} + \tau v e^{\tau v}}
                   = -\frac{e^{\tau\varepsilon\lambda^\alpha}(1 + \tau \varepsilon \alpha \lambda^\alpha)}{e^{\tau v} - \tau\lambda
                      e^{\tau \varepsilon \lambda^\alpha}} \\
                   &= -\frac{1 + \tau \varepsilon \alpha \lambda^\alpha}{e^{\tau v} e^{-\tau \varepsilon \lambda^\alpha} - \tau\lambda}
                   = \frac{1 + \tau \varepsilon \alpha \lambda^\alpha}{\frac{\lambda}{v} + \tau\lambda}
                   = \frac{\frac{1}{\lambda} + \tau \varepsilon \alpha \lambda^{\alpha - 1}}{\frac{1}{v} + \tau}.
            \end{split}
        \end{equation}
        Hence, since $\alpha < 1$, and since $v \to \infty$, we conclude
        \begin{equation}
            v' \to 0 \qquad \mbox{ as } \lambda \to \infty.
            \notag
        \end{equation}
        This completes the proof of (\ref{eq3.7}) and thus the proof of the Theorem \ref{THEOREM_ILL_POSEDNESS_LOW_ORDER_REG_HEAT_EQUATION}.
    \end{proof}

    We can extend the ill-posedness result to some higher-order
    equations of the type
    \begin{equation} \label{eqn-higher}
    \partial_t^m u(t) = {\cal A} u(t-\tau) - \varepsilon (-{\cal A})^\alpha
    u(t),
    \end{equation}
    where $m\geq 2$. Making a similar ansatz as in the proof of
    Theorem \ref{theo-m1}, the corresponding equation for $\omega$ is
    given by
     \begin{equation}\label{eq17}
        \omega^m + \varepsilon \lambda^\alpha = - \lambda e^{- \tau \omega}
    \end{equation}
    Ansatz:
    \begin{equation}\label{eq18}
        \omega = y_1 + i y_2 = r e^{i \varphi_m}
    \end{equation}
    with
    \begin{equation}\label{eq19}
        \varphi_m := \frac{\pi}{2 m}
    \end{equation}
    fixed.

    Then
    \begin{equation}\label{eq20}
        \omega^m = i r^m,
    \end{equation}
    and
    (\ref{eq17}) turns into
    \begin{align}
        \label{eq21}
        \varepsilon \lambda ^\alpha &= -\lambda e^{-\tau y_1} \cos (\tau y_2), \\
        \label{eq22}
        \left(y^2_1 + y^2_2\right)^{m/2} &= \lambda e ^{-\tau y_1} \sin (\tau y_2).
    \end{align}
    Observing
    \begin{equation}\label{eq23}
        y_2 = \beta y_1,
    \end{equation}
    with
    \begin{equation}\label{eq24}
        \beta = \tan(\varphi_m) > 0,
    \end{equation}
    Equations (\ref{eq21}), (\ref{eq22}) turn into
    \begin{align}
        \label{eq25}
        \varepsilon \lambda ^\alpha &= - \lambda e^{-\tau x} \cos (\tau \beta x), \\
        (1 + \beta^2)^m x^{\tau m} &= \lambda e^{-\tau x} \sin (\tau \beta x),
        \label{eq26}
    \end{align}
    where
    \begin{equation}\label{eq27}
        x := y_1, \quad (y_2 = \beta x).
    \end{equation}
    Here, the condition $\alpha < 1$ is important to give (\ref{eq25}) sense as $\lambda \to \infty$.

    From (\ref{eq26}) we have
    \begin{equation}\label{eq28}
        0 < \lambda = \frac{(1 + \beta^2)^m x^{2 m} \, e^{\tau x}}{\sin (\tau \beta x)},
    \end{equation}
    if
    \begin{equation}\label{eq29}
        \sin (\tau \beta x) > 0.
    \end{equation}
    Plugging (\ref{eq28}) into (\ref{eq25}) yields
    \begin{equation}\label{eq30}
        f_1(x):=  \varepsilon e^{\alpha \tau x} \big(\sin (\tau  \beta x) \big)^{1 - \alpha}
        = (1 + \beta^2)^{m (1 - \alpha)} x^{2 m (1 - \alpha)} \big(- \cos (\tau \beta x)\big)
        =:f_2(x),
    \end{equation}
    being well-defined if
    \begin{equation}\label{eq31}
        \cos(\tau \beta x) < 0.
    \end{equation}
    The equation (\ref{eq30}) has infinitely many solutions $x_k, k \in {\Bbb N}$,
    one in each interval
    \begin{equation}\label{eq32}
        I_k := \big(\tfrac{\pi + 4k \pi}{2 \tau \beta},
        \tfrac{\pi + 2 k\pi}{\tau \beta}  \big) \equiv (a_k,b_k)
    \end{equation}
    since $g := f_1 - f_2$ satisfies
    \begin{equation}\label{eq33}
        g(a_k)=f_1(a_k) >0 > -f_2(b_k) = g(b_k).
    \end{equation}
    Hence
    \begin{equation}
        {\rm Re}\, \omega_k = x_k \rightarrow \infty
        \notag
    \end{equation}
    and for $\lambda_k$, determined by (\ref{eq28}), we have
    \begin{equation}
        \lambda_k \rightarrow \infty.
        \notag
    \end{equation}
This way, the eigenvalues are not arbitrary, but we can define in what follows an
associated operator $\mathcal{A}$, for which we then have
the ill-posedness result related to equation (\ref{eqn-higher}). The
desired operator ${\cal A}$ can be chosen as
    \begin{equation}
        {\cal A}: D({\cal A}) \subset {\cal H} \rightarrow {\cal H} \notag
    \end{equation}
    in a Hilbert space ${\cal H}$ with a complete orthonormal system $(\Phi_k)_k \subset {\cal H} $ satisfying
    \begin{equation}
        \begin{split}
            {\cal A} \Phi_k &:= (-\lambda_k) \Phi_k, \\
            D({\cal A}) &= \Big\{u \in {\cal H} \Big| \sum^\infty_{k =1} \lambda^2_k |\langle u, \Phi_k \rangle|^2 < \infty\Big\}, \\
            {\cal A} u &= \sum^\infty_{k = 1} (-\lambda_k) \langle u, \Phi_k \rangle \Phi_k.
        \end{split}
        \notag
    \end{equation}
    For the special case $m = 2$ we can prove a similar result as for the case $m = 1$,
    i.e., we may prescribe the sequence $(-\lambda_n)_n$ of eigenvalues.
    Without loss of generality, we may assume $\varepsilon = \tau = 1$.
    The characteristic relation
    \begin{equation}
        \omega^2 + \lambda^\alpha = - \lambda e^{- \omega}
        \notag
    \end{equation}
    is, for
    \begin{equation}
        \omega = y_1 + i y_2, \quad (y_j \in {\Bbb R}),
        \notag
    \end{equation}
    equivalent to
    \begin{equation}\label{eq34}
        y^2_1 - y^2_2 + \lambda^\alpha = - \lambda e^{-y_1} \cos (y_2),
    \end{equation}
    \begin{equation}\label{eq35}
        2 y_1 y_2 = \lambda e^{-y_1} \sin (y_2).
    \end{equation}
    Looking for solutions satisfying
    \begin{equation}\label{eq36}
        y_2 \in [\pi/2, \pi),
    \end{equation}
    Equation (\ref{eq35}) is equivalent to
    \begin{equation}\label{eq37}
        y_1 e^{y_1} \; \tfrac{2 y_2}{\sin (y_2)} = \lambda.
    \end{equation}
    Defining
    \begin{equation}
        \widetilde{h}_2 : [\pi/2, \pi) \rightarrow [\pi, \infty), \quad \widetilde{h}_1 : [0, \infty) \rightarrow [0, \infty)
        \notag
    \end{equation}
    via
    \begin{equation}
        \widetilde{h}_1 (y_1) := y_1 e^{y_1}, \quad \widetilde{h}_2 (y_2) := \tfrac{2 y_2}{\sin (y_2)},
        \notag
    \end{equation}
    we have that $\widetilde{h}^\prime_1 > 0$ and $\widetilde{h}^\prime_2 > 0$.
    Let
    \begin{equation}
        h_1 := \widetilde{h}^{- 1}_1 : [0, \infty) \rightarrow [0, \infty),
        \quad h_2 := \widetilde{h}^{- 1}_2 : [\pi, \infty) \rightarrow [\tfrac{\pi}{2}, \pi)
    \end{equation}
    satisfy
    \begin{equation}\label{eq38}
        h_1 (0) = 0,\; \lim_{z \to \infty} h_1 (z) =  \infty, \quad h_2 (\pi) = \tfrac{\pi}{2}, \;
        \lim_{\eta \to \infty} h_2(\eta) = \pi.
    \end{equation}
    According to (\ref{eq37}), one has to fulfill
    \begin{equation}
        \widetilde{h}_1 (y_1) \widetilde{h}_2 (y_2) = \lambda,
        \notag
    \end{equation}
    hence, allowing
    \begin{equation}
        \pi \leq \widetilde{h}_2 (y_2) < \infty,
        \notag
    \end{equation}
    one requires
    \begin{equation}
        0 < \widetilde{h}_1 (y_1) \leq \tfrac{\lambda}{\pi}.
        \notag
    \end{equation}
    Therefore, $h_1$ is considered restricted to
    \begin{equation}\label{eq39}
        h_1 \colon (0, \tfrac{\lambda}{\pi}] \rightarrow (0, h_1 (\tfrac{\lambda}{\pi})].
    \end{equation}
    Denoting
    \begin{equation}
        y_1 = h_1 (z), \quad y_2 = h_2  (\eta),
        \notag
    \end{equation}
    Equations (\ref{eq34}), (\ref{eq35}) turn into
    \begin{equation}\label{eq40}
        h^2_1 (z) - h^2_2 (z) + \lambda^\alpha = -\lambda e^{-h_1 (z)} \cos (h_2 (\eta)),
    \end{equation}
    \begin{equation}\label{eq41}
        z \,\eta = \lambda,
    \end{equation}
    for
    \begin{equation}\label{eq42}
        (z, \eta) \in G_\lambda := (0, \tfrac{\lambda}{\pi}] \times [\pi, \infty).
    \end{equation}
    We look for solutions $(z_n, \eta_n) \in G_{\lambda_n}$ to (\ref{eq40}), (\ref{eq41})
    satisfying $z_n \to \infty$ as $n \to\infty$.

    That is, using (\ref{eq41}), we wish to solve
    \begin{equation}\label{eq43}
        {\cal F} (\lambda_n, z) := h^2_1 (z) - h^2_2 (\tfrac{\lambda}{z}) + \lambda^\alpha + \lambda e^{-h_1 (z)} \cos \big(h_2 (\tfrac{\lambda}{z}) \big) = 0.
    \end{equation}
    Since $\lim\limits_{\eta \to \infty}^{} h_2 (\eta) = \pi$, we have
    \begin{equation}\label{eq44}
        \lim_{z \to 0} {\cal F} (\lambda_n, z) = - \pi^2 + \lambda^\alpha_n - \lambda_n < 0
    \end{equation}
    if $n$ is large enough, $n \geq n_0$ for some $n_0 \in {\Bbb N}$.

    Since $h_2 (\pi) = \tfrac{\pi}{2}$, we have
    \begin{equation}
        \label{eq45}
        \begin{split}
            {\cal F} (\lambda_n, \tfrac{\lambda_n}{\pi}) &= h^2_1 (\tfrac{\lambda_n}{\pi}) - \tfrac{\pi^2}{4} + \lambda^\alpha_n \\
            &\geq h^2_1 (\tfrac{\lambda_n}{\pi}) - \tfrac{\pi^2}{4} > 0
        \end{split}
    \end{equation}
    if $n \geq n_1$ for some $n_1 \in {\Bbb N}$.

    On the strength of (\ref{eq43}), (\ref{eq44}), we conclude that, if $n \geq n^{\ast} := \max\{n_0, n_1\}$,
    \begin{equation}\label{eq46}
        \exists z \equiv z_n \equiv z (\lambda_n) \in (0, \tfrac{\lambda_n}{\pi}) :
         {\cal F} (\lambda_n, z_n) = 0.
    \end{equation}
    There remains to prove that there exists (at least a) subsequence
    $(\hat{z}_n)$ of $(z_n)_n$ such that $\hat{z}_n \rightarrow \infty$.
    For this purpose  we observe from (\ref{eq45}) that
    \begin{equation}\label{eq47}
        \frac{h^2_1 (z_n)}{\lambda_n} - \frac{h_2(\tfrac{\lambda}{z_n})}{\lambda_n} +
        \frac{\lambda_n^\alpha}{\lambda} + e^{-h_1 (z_n)}
         \cos \big(h_2 (\tfrac{\lambda_n}{z_n}) \big) = 0.
    \end{equation}
    Assuming
    \begin{equation}\label{eq48}
        \sup_{n \in\, {\Bbb N}} h_1 (z_n) < \infty,
    \end{equation}
    we conclude from (\ref{eq46}), using the boundedness of $h_2$ and, in particular, $\alpha < 1$,
    \begin{equation}\label{eq49}
         \lim_{n \to \infty} e^{-h_1 (z_n)} \cos \big(h_2 (\tfrac{\lambda_n}{z_n})\big) = 0,
    \end{equation}
    implying, by assumption (\ref{eq46}),
    \begin{equation}
        \lim_{n \to \infty} \cos h_2 (\tfrac{\lambda_n}{z_n}) = 0, \notag
    \end{equation}
    or,
    \begin{equation}\label{eq50}
        \frac{\lambda_n}{z_n} \rightarrow \pi
    \end{equation}
    implying that $(z_n)_n$ is unbounded, which implies, for a subsequence $(\hat{z}_n)_n$,
    that $\hat{z}_n \to \infty$ which is a contradiction to the assumption (\ref{eq46}).
    Therefore $(h_1 (z n))_n$  is unbounded, implying the existence of a
    subsequence $(\hat{z}_n)_n$ with $\hat{z}_n \to \infty$.\\
    Thus, we have proved
    \begin{theorem} \label{theo-mgeq2}
	\begin{itemize}
		\item[\textit{(i)}] For $m\geq 2$ there are operators $\cal{A}$ associated to (\ref{eqn-higher}) for which the problem is ill-posed if $\alpha < 1$.
		\item[\textit{(ii)}] For $m=2$, a result corresponding to Theorem \ref{theo-m1} holds true.
	\end{itemize}
    \end{theorem}
    \begin{remark}
        The arguments do not carry over to the case $m = 1$ since (\ref{eq47}) does no longer
        follow (instead $e^{-h_1 (z_{n})} \cos h_2 (\frac{\lambda_n}{z_{n}}) \rightarrow -1$).
    \end{remark}
    \begin{remark}
        The case $m \geq 3$ and {\em prescribing} $(-\lambda_n)_n$ is still open.
    \end{remark}

\section{Physical Example}
    In this last section, we apply the explicit solution representation from Section \ref{SUBSECTION_REPRESENTATION_OF_SOLUTIONS}
    to solve a physical problem arising from microscale thermal transport phenomena in thin metal films.
    A kinetic description of the latter can be derived from the Boltzmann equation for electrons and phonons
    (see \cite{AnKaPe1974}, \cite{KoNieJaSchKrPo1998}, \cite{TaZh1998}).

    We consider a 50 nm thin gold film occupying the interval $(0, L)$ of the real line (i.e., $L = 50 \cdot 10^{-9}$ [m]).
    Let $u, \theta, q$ denote the electron energy density, electron temperature, and electron heat flux, respectively.
    For simplicity, we assume the phonon temperature $\theta_{l} \equiv 300$ [K] to be constant
    and the phonon heat flux $q_{l} \equiv 0$ to vanish.
    If the electron gas is in equilibrium, its energy density is related to the electron temperature as
    $u = \tfrac{\gamma}{2} \theta^{2}$ for a positive $\gamma$.
    Performing linearization around $\theta_{l}$, we obtain the following constitutive equation:
    \begin{equation}
        u = c_{e} \theta \notag
    \end{equation}
    where $c_{e} := \gamma \theta_{l}$ is the electron heat capacity.
    The film is assumed to undergo a short pump laser pulse applied to its left surface (i.e., $x = 0$)
    causing an increase in electron temperature (see \cite{KoNieJaSchKrPo1998}, \cite{MaYi2012}, \cite{QuTie1992}).
    The absorption of the laser radiation is modeled by a source term $f$
    (cf. \cite{MaYi2012}, \cite{QuTie1992}). See Table \ref{TABLE_LASER_CONSTANTS} for details.
    \begin{table}[h!]
        \centering
        \begin{tabular}{cccc}
            Notation & Units & Value & Description \\
            \hline\hline
            $\gamma$ & $\mathrm{J}\,\mathrm{m}^{-3}\,\mathrm{K}^{-2}$ & $67.6 \cdot 10^{-3} $ & electron heat capacity increase per degree ${}^{\circ}\mathrm{K}$ \\
            $c_{e}$ & $\mathrm{J}\,\mathrm{m}^{-3}\,\mathrm{K}^{-1}$ & $2.1 \cdot 10^{4}$ & electron heat capacity \\
            $\tau$ & $\mathrm{s}$ & $26 \cdot 10^{-15}$ & electron relaxation time \\
            $\lambda$ & $\mathrm{W}\,\mathrm{m}^{-1}\,\mathrm{K}^{-1}$ & 315 & electron thermal conductivity \\
            $G$ & $\mathrm{W}\,\mathrm{m}^{-3}\,\mathrm{K}^{-1}$ & $2.6 \cdot 10^{16}$ & electron-lattice coupling constant
        \end{tabular}
        \caption{Material properties of Au (gold) \label{TABLE_CONSTANTS}}
    \end{table}
    Replacing Cattaneo's law with a regularized delay law and neglecting the equations for the phonon variables,
    Equation (53) from \cite{TaZh1998} is reduced to
    \begin{equation}
        \begin{split}
            c_{e} \partial_{t} \theta(t, x) + \partial_{x} q(t, x) + G \cdot (\theta(t, x) - \theta_{l}) &= f(t, x) \text{ for } (t, x) \in (0, T) \times (0, L), \\
            q(t, x) + \varepsilon \lambda \partial_{x} \theta(t, x) + \lambda \partial_{x} \theta(t - \tau, x) &= 0 \text{ for } (t, x) \in (0, T) \times (0, L)
        \end{split}
        \label{EQUATION_ELECTRON_HEAT_CONDUCTION_SYSTEM}
    \end{equation}
    where $\rho$ is the density, $c_{\rho}$ the electron heat capacity, $G$ electron-lattice coupling factor,
    and $\lambda$ electron thermal conductivity.
    Eliminating $q$ from Equation (\ref{EQUATION_ELECTRON_HEAT_CONDUCTION_SYSTEM}) yields
    \begin{equation}
	\begin{split}
		\partial_{t} \theta(t, x) - \tfrac{\varepsilon \lambda}{c_{e}} \partial_{xx} \theta(t, x)
        	+ \tfrac{G}{c_{e}} \theta(t, x) - \tfrac{\lambda}{c_{e}} \partial_{xx} \theta(t - \tau, x)
        	= \tfrac{1}{c_{e}} f(t, x) + \tfrac{G}{c_{e}} \theta_{l} \\
		\hfill \text{ for } (t, x) \in (0, T) \times (0, L).
	\end{split}
        \label{EQUATION_ELECTRON_HEAT_CONDUCTION_SYSTEM_SIMPLIFIED}
    \end{equation}
    To close the equation, we prescribe homogeneous Neumann boundary conditions
    \begin{equation}
        \partial_{x} \theta(t, 0) = \partial_{x} \theta(t, L) = 0 \text{ for } (t, x) \in (0, T) \times (0, L)
        \label{EQUATION_ELECTRON_HEAT_CONDUCTION_SYSTEM_SIMPLIFIED_BC}
    \end{equation}
    modeling the insulation of film surface
    and the initial condition
    \begin{equation}
        \theta(t, x) \equiv \theta^{0} \text{ for } (t, x) \in (-\tau, 0) \times (0, L), \quad
        \theta(0, x) \equiv \theta^{0} \text{ for } x \in (0, L)
        \label{EQUATION_ELECTRON_HEAT_CONDUCTION_SYSTEM_SIMPLIFIED_IC}
    \end{equation}
    with $\theta^{0} \equiv 300$ [$\mathrm{K}$].
    \begin{table}[h!]
        \centering
        \begin{tabular}{cccc}
            Notation & Units & Value & Description \\
            \hline\hline
            $r_{f}$ & --- & 0.94 & reflectivity \\
            $t_{p}$ & $\mathrm{s}$ & $96 \cdot 10^{-15}$ & laser peak time \\
            $\alpha^{-1}$ & $\mathrm{m}$ & $15 \cdot 10^{-9}$ & laser radiation penetration depth \\
            $J$ & $\mathrm{J}\, \mathrm{m}^{-2}$ & 150 & total laser energy over the spot cross-section \\
        \end{tabular}
        \caption{Laser source term $f(t, x) = 0.94 \frac{1 - r_{f}}{t_{p}} \alpha J \exp\left(-x\alpha - 2.77 \big(\tfrac{t}{t_{p}}\big)^{2}\right)$
        \label{TABLE_LASER_CONSTANTS}}
    \end{table}

\begin{figure}[h!]
        \centering
        \includegraphics[scale = 0.35]{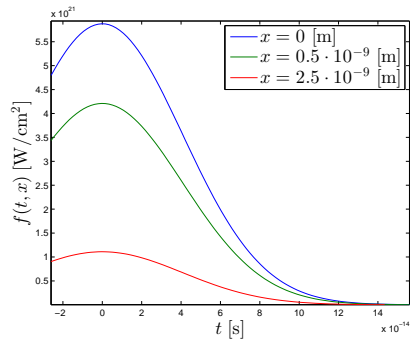}
        \caption{Laser intensity at different space points \label{FIGURE_LASER_INTENSITY}}
    \end{figure}

    Figure \ref{FIGURE_LASER_INTENSITY} displays the laser intensity at three different points in the film.

\begin{figure}[h!]
        \centering
        \begin{tabular}{cc}
            \includegraphics[scale = 0.35]{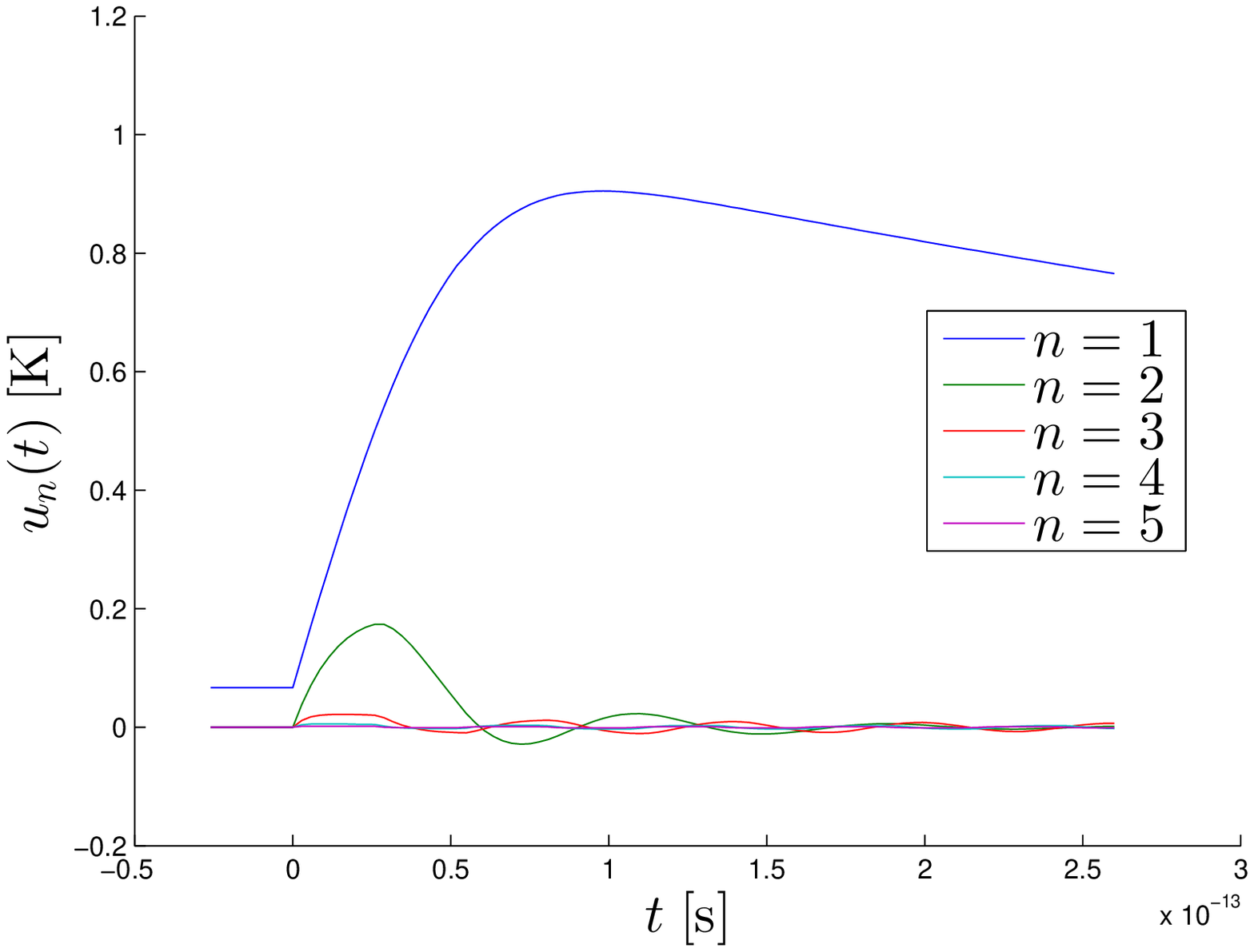} & \includegraphics[scale = 0.35]{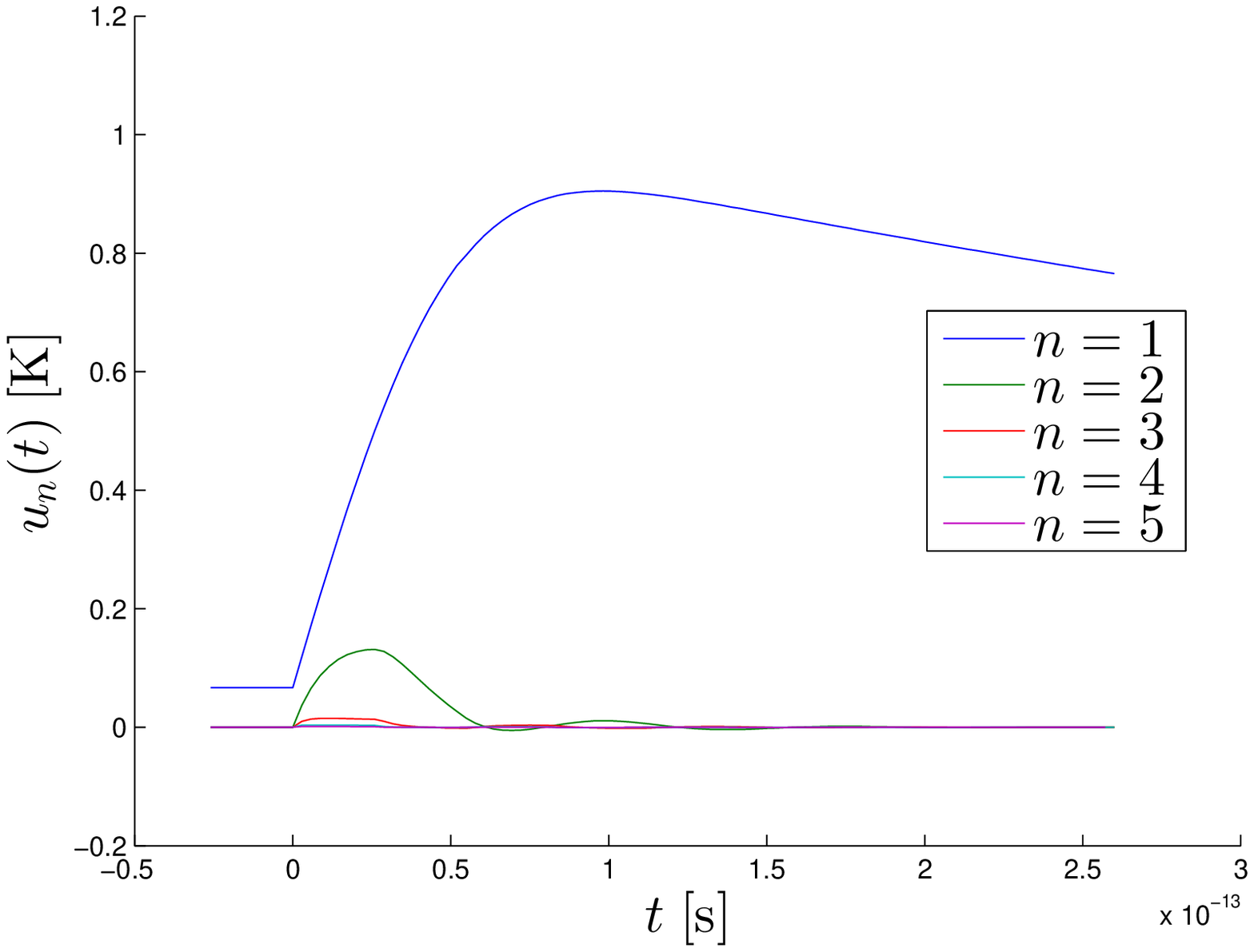} \\
            $\varepsilon = 1$ & $\varepsilon = 1.5$
        \end{tabular}
        \caption{Time-dependent Fourier coefficients $u_{n}$ \label{FIGURE_SOL_N}}
    \end{figure}

    Our theory from Section \ref{SECTION_REGULARIZED_HEAT_CONDUCTION_WITH_DELAY}
    does not directly apply to the problem
    (\ref{EQUATION_ELECTRON_HEAT_CONDUCTION_SYSTEM_SIMPLIFIED})--(\ref{EQUATION_ELECTRON_HEAT_CONDUCTION_SYSTEM_SIMPLIFIED_IC})
    since Neumann and not Dirichlet boundary conditions are prescribed.
    Nonetheless, an analogous explicit representation formula as the one obtained in Section \ref{SUBSECTION_REPRESENTATION_OF_SOLUTIONS}
    directly applies to this new problem
    with $(\lambda_{n})_{n}$ and $(\phi_{n})_{n}$ replaced
    by the eigenvalues and orthonormal eigenfunctions of Neumann-Laplacian on $\Omega := (0, L)$.
    The latter read as
    \begin{equation}
        \lambda_{n} = \tfrac{\pi^{2} (n - 1)^{2}}{L^{2}}, \quad
        \phi_{n}(x) =
        \left\{\begin{array}{cl}
            \tfrac{1}{\sqrt{L}}, & n = 1, \\
            \sqrt{\tfrac{2}{L}} \cos\big(\tfrac{(n - 1) \pi x}{L}\big), & n > 1,
        \end{array}\right.
        \quad x \in [0, L], \quad\quad \text{ for } n \in \NN.
        \notag
    \end{equation}
    Using $f$ given in Table \ref{TABLE_LASER_CONSTANTS}, we compute
    \begin{equation}
        \begin{split}
            \langle f(t, \cdot), \phi_{n}\rangle \approx
            \left\{\begin{array}{cl}
                \tfrac{0.94 J}{t_{p} \sqrt{L}} \exp\Big(\tfrac{-\alpha t_p^2 L + 277 t^2}{t_p^2}\Big) (-1 + r_f) (e^{\alpha L} - 1), & n = 1, \\
                \tfrac{1.33 \alpha^2 L^{3/2}}{t_{p} \alpha^2 L^2 + 9.87 n^2} \exp\Big(\tfrac{-\alpha t_p^2 L + 277 t^2}{t_p^2}\Big) (-1 + r_f) (e^{\alpha L} + (-1)^{n}), & n > 1.
            \end{array}\right.
            \notag
        \end{split}
    \end{equation}
    Plugging these data into Equation (\ref{EQUATION_DELAY_ODE_N})
    and using the solution formula (\ref{EQUATION_DELAY_ODE_N_SOL_FORMULA}),
    we can explicitly compute $(u_{n})_{n}$.
    We performed this
    using Simpson's quadrature formula to numerically evaluate the integrals.
    The solutions are plotted in Figure \ref{FIGURE_SOL_N}.

\begin{figure}[h!]
        \centering
        \begin{tabular}{cc}
            \includegraphics[scale = 0.35]{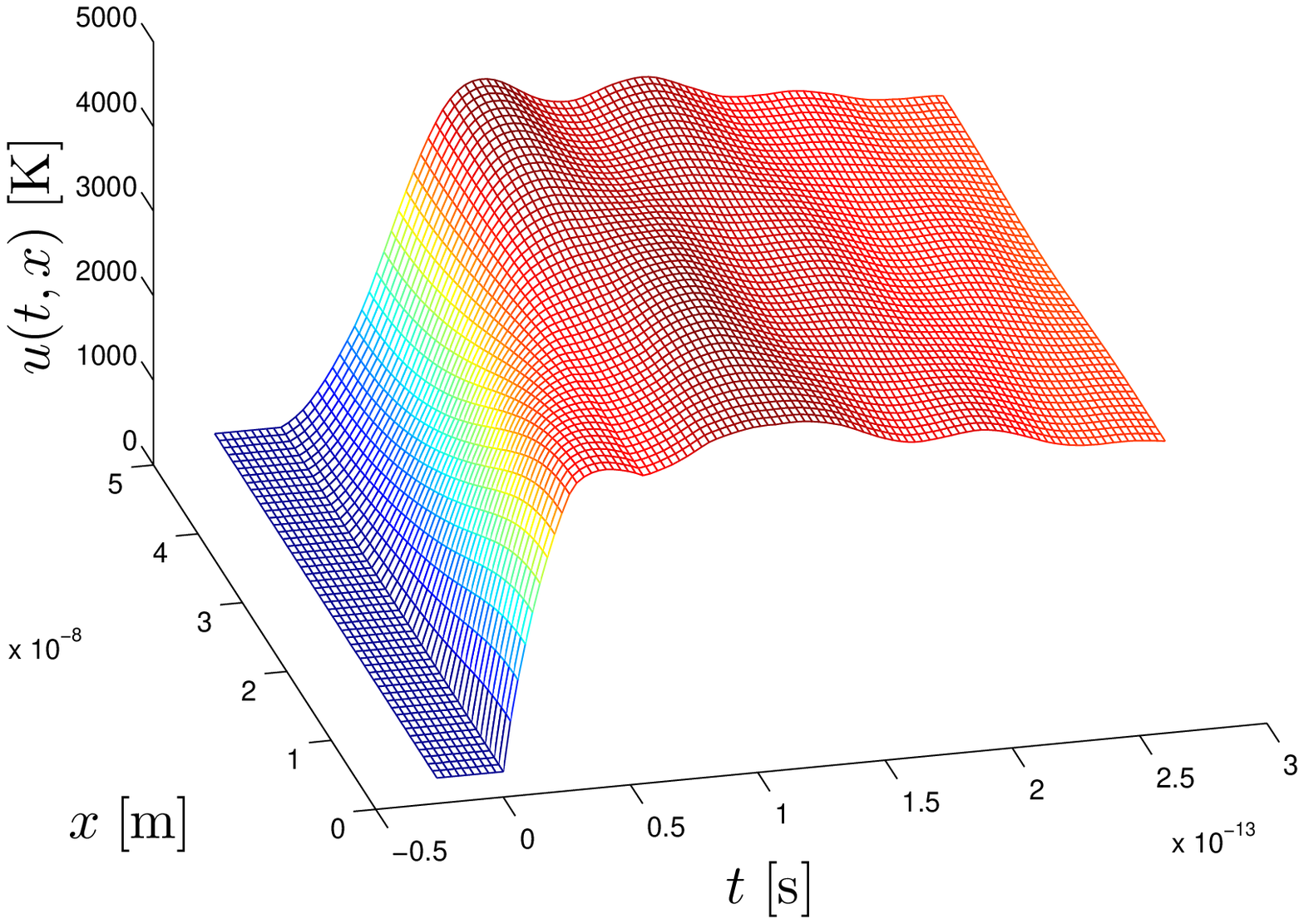} & \includegraphics[scale = 0.35]{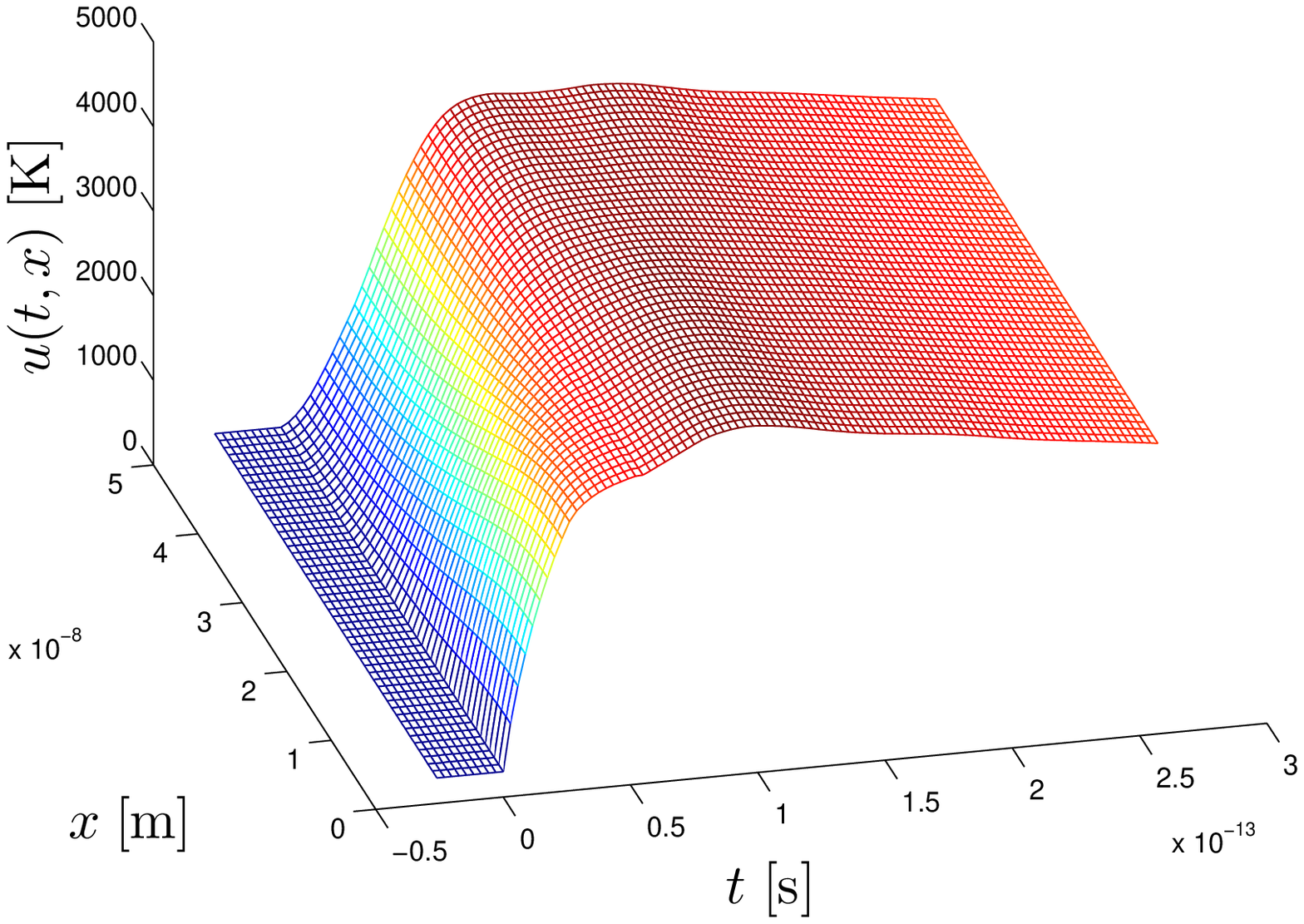} \\
            $\varepsilon = 1$ & $\varepsilon = 1.5$
        \end{tabular}
        \caption{Numerical solution \label{FIGURE_NUMERICAL_SOLUTION}}
    \end{figure}

    Plugging this numerical solution into Equation (\ref{EQUATION_DELAY_ODE_SOL_FORMULA})
    and considering first $n \leq 5$ terms in the series,
    we finally obtain a numerical solution plotted in Figure \ref{FIGURE_NUMERICAL_SOLUTION}.
    Note that these first five terms provide a very accurate approximation since higher Fourier coefficients practically vanish.
    The solution function has a peak somewhere at $\hat{t}_{p} = (95 \pm 5) \cdot 10^{-15} [\mathrm{s}]$
    which is close to the expected peak value $t_{p} = 96 \cdot 10^{-15}$.

    When $\varepsilon$ increases, the solution function becomes smoother.
    For $\varepsilon < 1$, e.g., $\varepsilon = 0.5$, the solution function becomes very rough
    due to the high volatility of Fourier coefficients.
    This observation suggests that the regularization parameter $\varepsilon$ should be selected
    to achieve best fit with experimental measurements.

\section{Appendix: Semi-discrete Lebesgue and Sobolev Spaces}
    Let $X$ be a Hilbert space
    and let $a = i\tau$, $b = j\tau$ for some $\tau > 0$ and $i, j \in \NN_{0}$ with $i < j$.
    We introduce the following semi-discrete Hilbert spaces
    \begin{equation}
		\begin{split}
			L^{2}_{\tau}\big((a, b), X\big) &:= \big\{u = (u_{i}, \dots, u_{j}) \,\big|\, u_{k} \in L^{2}\big((0, \tau), X\big) \text{ for } i \leq k \leq j\big\}, \\
			H^{1}_{\tau}\big((a, b), X\big) &:= \big\{u = (u_{i}, \dots, u_{j}) \,\big|\, u_{k} \in H^{1}\big((0, \tau), X\big), u_{k}(\tau) = u_{k+1}(0) \\
			&\hspace{5.5cm} \text{ for } i \leq k < k + 1 \leq j\big\}
		\end{split} \notag
    \end{equation}
    endowed with the standard product topology, i.e.,
    \begin{equation}
        \|u\|_{L^{2}_{\tau}((a, b), X)}^{2} = \sum_{k = i}^{j} \|u_{k}\|_{L^{2}((0, \tau), X)}^{2}, \quad
        \|u\|_{H^{1}_{\tau}((a, b), X)}^{2} = \sum_{k = i}^{j} \|u_{k}\|_{H^{1}((0, \tau), X)}^{2}. \notag
    \end{equation}
    Note that due to the continuity of the embedding
    $H^{1}\big((0, \tau), X\big) \hookrightarrow \mathcal{C}^{0}\big([0, \tau], X\big)$
    the space $H^{1}_{\tau}\big((a, b), X\big)$ is well-defined.

    Next, we consider the mapping
    \begin{equation}
        R \colon L^{2}\big((a, b), X\big) \to L^{2}_{\tau}\big((a, b), X\big), \quad
        u \mapsto (r_{i} u, \dots, r_{j} u) \notag
    \end{equation}
    with $(r_{k} u)(s) = u((k - 1) \tau + s)$
    for $s \in [0, \tau]$, $k = i, \dots, j$.
    Obviously, $R$ is an isomorphism.
    Moreover, the following assertion holds true.
    \begin{lemma}
        \label{LEMMA_SEMIDISCRETE_CONTINUOUS_EQUIVALENCE}
        $u \in H^{1}\big((a, b), X\big)$ is true if and only if
        $Ru \in H^{1}_{\tau}\big((a, b), X\big)$.
    \end{lemma}
    \begin{proof}
        If $b - a = \tau$, the claim is trivial.
        Due to Sobolev embedding, the implication
        $u \in H^{1}\big((a, b), X\big) \Rightarrow Ru \in H^{1}_{\tau}\big((a, b), X\big)$
        also trivially follows.
        To show the converse, due to \cite{Ad1975},
        it suffices to prove that
        \begin{equation}
            Ru \in H^{1}_{\tau}\big((a, b), X\big) \Rightarrow
            u \in W^{1, 2}\big((a, b), X\big). \notag
        \end{equation}
        For $\phi \in \mathcal{C}^{0}\big((a, b), X\big)$, we obtain
        \begin{equation}
            \begin{split}
                \int_{a}^{b} u(t) \cdot \partial_{t} \phi(t) \mathrm{d}t &=
                \sum_{k = i}^{j-1} \int_{i\tau}^{(i+1) \tau} u_{k}(t) \cdot \partial_{t} \phi(t) \mathrm{d}t \\
		& =-\sum_{k = i}^{j-1} \int_{i\tau}^{(i+1) \tau} \partial_{t} u_{k}(t) \cdot \phi(t) \mathrm{d}t +
                u(t) \phi(t)|_{t = k\tau}^{t = (k+1) \tau} \\
                &= -\int_{a}^{b} \Big(\sum_{k = i}^{j-1} \partial_{t} u_{k}(t) \chi_{(k\tau, (k+1) \tau)}(t)\Big) \cdot \phi(t) \mathrm{d}t
                =: -\int_{a}^{b} \partial_{t} u(t) \cdot \phi(t) \mathrm{d}t \notag
            \end{split}
        \end{equation}
        since the boundary terms vanish
        due to the definition of $H^{1}_{\tau}\big((a, b), X\big)$ and the fact that $\phi(a) = \phi(b) = 0$.
        Finally, we observe
        \begin{equation}
            \int_{a}^{b} \|\partial_{t} u(t)\|_{X}^{2} \mathrm{d}t =
            \sum_{k = i}^{j-1} \int_{i\tau}^{(i+1) \tau} \|\partial_{t} u(t)\|_{X}^{2} \mathrm{d}t < \infty. \notag
        \end{equation}
        Thus, $u \in W^{1,2}\big((a, b), X\big) = H^{1}\big((a, b), X\big)$.
    \end{proof}

\end{document}